\input amstex
\loadbold
\documentstyle{amsppt}
\magnification=1200

\input xy
\xyoption{all}

\pagewidth{5.65truein}
\pageheight{9.2truein}
\hcorrection{0.8truecm}

\def\nmb#1#2{#2}         % used for renumbering, TeX should ignore.
\def\cit#1#2{\ifx#1!\cite{#2}\else#2\fi} %for citing references
\def\totoc{}             %= to table of content, invoked by kms-book.sty
               % for producing index, invoked by kms-book.sty
\def\ign#1{}             %=ignore, invisible entry for the index only

\BlackBoxes

\define\pr{0}

\define\sm{\operatorname{sm}}
\define\rat{\operatorname{rat}}
\define\dom{\operatorname{dom}}

\define\Div{\operatorname{Div}}
\define\dvs{\operatorname{div}}
\define\supp{\operatorname{supp}}
\define\ve{\varepsilon}
\define\vph{\varphi}

\define\la{\lambda}

\define\ph{\varphi}

\define\Ga{\Gamma}

\redefine\i{^{-1}}
\define\row#1#2#3{#1_{#2},\ldots,#1_{#3}}
\define\x{\times}

\def\today{\ifcase\month\or
 January\or February\or March\or April\or May\or June\or
 July\or August\or September\or October\or November\or December\fi
 \space\number\day, \number\year}

\topmatter
\title \vskip 5mm Invariant tensor fields and
orbit varieties\\ for
       finite
       algebraic transformation groups
\endtitle
\author Mark Losik, Peter W. Michor, Vladimir L. Popov
\endauthor
\rightheadtext{Invariant tensor fields for
algebraic finite group actions}
\leftheadtext{M.\ Losik, P. W.\ Michor, V.
L. Popov}
\address
M. Losik: Saratov State University, ul. Astrakhanskaya, 83,
410026 Saratov, Russia
\endaddress
\email LosikMV\@info.sgu.ru \endemail
\address
P\. W\. Michor: Institut f\"ur Mathematik, Universit\"at Wien,
Strudlhofgasse 4, A-1090 Wien, Austria; {\it and:}
Erwin Schr\"odinger Institut f\"ur Mathematische Physik,
Boltzmanngasse 9, A-1090 Wien, Austria
\endaddress
\email Peter.Michor\@esi.ac.at \endemail
\address
V.L\. Popov: Steklov Mathematical Institute,
Russian Academy of Sciences, Gub\-ki\-na 8,
Moscow, 117966, Russia
\endaddress
\email popov\@ppc.msk.ru \endemail

\dedicatory{
To C. S. Seshadri on
the occasion of his 70th birthday}
\enddedicatory

\thanks{M.L. and P.W.M. were supported
     by ``Fonds zur F\"orderung der
     wissenschaftlichen For\-schung,
     Projekt P~14195~MAT".}
\endthanks
 \abstract
Let $X$ be a smooth algebraic variety endowed
with an action of  a finite group $G$ such that
there exists a geometric quotient
$\pi_X:X\rightarrow X/G$. We characterize
rational tensor fields $\tau$ on $X/G$ such that
the {\it pull back} of $\tau $ is regular on $X$:
these are precisely all $\tau$ such that
$\dvs_{R_{X\!/\!G}}(\tau)\geqslant 0$ where
$R_{X\!/\!G}$ is the {\it reflection divisor} of
$X/G$ and $\dvs_{R_{X\!/\!G}}(\tau)$ is the {\it
$R_{X\!/\!G}$-divisor} of $\tau$. We give some
applications, in particular to a generalization
of Solomon's theorem. In the last section we show
that if $V$ is a finite dimensional vector space
and $G$ a finite subgroup of
$\operatorname{GL}(V)$, then each automorphism
$\psi$ of $V/G$ admits a biregular lift $\varphi:
V\rightarrow V$ provided that $\psi$ maps the
regular stratum to itself and
$\psi_*(R_{X\!/\!G})=R_{X\!/\!G}$.
\endabstract
\endtopmatter

\document

\head\totoc\nmb0{\bf 1}. {\bf
Introduction}
\endhead
Let $X$ be a smooth algebraic variety endowed
with an action of a finite group $G$. Assume that
there exists a geometric quotient $\pi^{}_{X}:
X\rightarrow X/G$ (this is always the case if $X$
is quasi-projective, cf. Subsection 2.4).

In this paper we study the
interrelations between rational tensor
fields on $X$ and $X/G$. If
$\tau$ is a rational tensor field on $X/G$,
we define a $G$-invariant
rational tensor field $\pi^{*}_{X}(\tau)$
on $X$ called the {\it pull back} of
$\tau$.
If $\theta$ is
a $G$-invariant rational tensor field
on $X$, we define  a
rational tensor field $\pi^{}_{X*}(\theta)$
on $X/G$ called the {\it push forward} of $\theta$.
We have
$\pi^{}_{X*}(\pi^{*}_{X}(\tau))=\tau$ and
$\pi^{*}_{X}(\pi^{}_{X*}(\theta))=\theta$.

Given this, we consider the following
problem: Let $\tau$ be a rational tensor
field on $X/G$. When is the pull back
$\pi^{*}_{X}(\tau)$ regular on $X$?

To that end we consider the Luna stratification
of $X/G$. Let $(X/G)_1$ be the union of all
codimension 1 strata. We show that $(X/G)_1$ is
contained in the smooth locus of $X/G$. Let
$(X/G)_1=(X/G)_1^1\cup\ldots\cup (X/G)_1^d$ be
the decomposition into irreducible components. We
show that for any $l=1,\ldots, d$ and $z\in
(X/G)^l_1$, $x\in \pi^{-1}_{X}(z)$ the stabilizer
of $x$ is a cyclic group whose image under the
slice representation is generated by a
pseudo-reflection of order $r^{}_l$ depending
only on $l$. We encode this information into the
{\it reflection divisor}
$R_{X/G}:=r_1R_1+\ldots+r_dR_d\in \Div(X/G)$
where $R_l$ is the prime divisor whose support is
the closure of $(X/G)^l_1$. Further, for any
nonzero rational tensor field $\tau $ on $X/G$
and positive divisor $B\in \Div(X/G)$ we define
$\dvs(\tau)\in \Div(X/G)$, the {\it divisor of}
$\tau$, and $\dvs_B(\tau)\in \Div(X/G)$, the
$B$-{\it divisor of} $\tau$. Locally $\dvs(\tau)$
is the minimum of divisors of the component
functions of $\tau$ in the framings of tangent
and cotangent bundles. The divisor $\dvs_B(\tau)$
is obtained from $\dvs(\tau)$ by means of some
``modification along" $B$ (see Subsection 3.6).

Our main result is that $\pi^*_X(\tau)$
is regular on $X$ if and only if
$\dvs_{R_{X\!/\!G}}(\tau)\geqslant 0$.

As a corollary we obtain a generalization
of Solomon's theorem \cite{So}, and prove
that the push forward of a $G$-invariant
regular tensor field on $X$ that is skew
symmetric with respect to the covariant
entries is regular on the smooth locus of
$X/G$.

Another application pertains to the case
where $X$ is a vector space $V$ with a
linear action of $G$: we obtain a
characterization of $G$-invariant
polynomials on $V^{\oplus s}$ in terms of
rational multi-symmetric covariant tensor
fields on $V/G$.

In the last section we prove that any
automorphism $\psi$ of the algebraic
variety $V/G$ such that $\psi
((V/G)_0)\subseteq (V/G)_0$ and
$\psi_*(R_{V\!/\!G})=R_{V\!/\!G}$ can be
lifted to an automorphism of the
algebraic variety $V$. The proof is based
on the relevant result in the analytic
setting, \cite{KLM}, so what we really
prove is that analytic lift is actually
algebraic.

Throughout in this paper we assume that
the base field $k$ is algebraically
closed of characteristic $0$. In the last
section we use the result from \cite{KLM}
that is proved for $k=\Bbb C$, so passing
to the general case is carried out by
Lefschetz's principle. However note that
our proofs of the results from Sections
2--4 could be extended {\it mutatis
mutandis} to the case when
$\operatorname{char} k$ is positive and
subject to some non-divisibility and
magnitude conditions; for this purpose
one should use the relevant replacement
of the slice theorem \cite{Lu} proved in
\cite{BR}.

This paper is the algebraic geometrical
companion of paper \cite{KLM} where
similar results were obtained for
analytic actions of finite groups.

\smallskip

We thank Yu.\,Neretin for helpful discussions.
Thanks are also due to the referee for
suggestions and comments which has led to
eliminating some inaccuracies and improvement in
the exposition.

\medskip

{\it Notation, terminology and
conventions}:

\smallskip
\roster \item"$\bullet\hskip 1.3mm     $" $\Cal
O_{x, X}$ and $T_{x, X}$ are respectively the
local ring and tangent space at a point $x$ of an
algebraic variety $X$. \item"$\bullet\hskip 1.3mm
$" $X_{\sm}$ is the smooth locus of $X$.
\item"$\bullet\hskip 1.3mm     $"
$\operatorname{Div}(X)$ is the Weil divisor group
of $X$. \item"$\bullet\hskip 1.3mm     $"
$\dvs(f)$ is the divisor of a rational function
$f$. \item"$\bullet\hskip 1.3mm $"$\supp(D)$ is
the support of $D\in \Div(X)$.
\item"$\bullet\hskip 1.3mm $"
  A positive divisor is called {\it prime} if it
  not sum of two positive divisors.
\item"$\bullet\hskip 1.3mm $" $m^{}_{C, D}$ is
the multiplicity of a prime divisor $C$ in a
divisor $D$. \item"$\bullet\hskip 1.3mm $"
$d_x\varphi$ is the differential of a morphism
$\varphi $ at a point $x$. \item"$\bullet\hskip
1.3mm     $" $G\cdot z$ and $G_z$ are
respectively the orbit and stabilizer of a point
$z$ of a set $Z$ endowed with an action of a
group $G$. \item"$\bullet\hskip 1.3mm     $"
$Z^S:=\{z\in Z\mid g\cdot z=z \text{ for all }
g\in S\}$. \item"$\bullet\hskip 1.3mm     $"
$|S|$ is the number of elements of a finite set
$S$.
\endroster

All group actions considered in this
paper are algebraic actions on algebraic
varieties.

\head\totoc\nmb0{\bf 2}. {\bf
Preliminaries}
\endhead

\subhead\nmb.{2.1.} Excellent morphisms
\endsubhead Recall that a morphism of
algebraic varieties $\ph:X\to Y$ is called {\it
\'etale at a point} $x\in X$ if the homomorphism
of local rings $\ph^*:\Cal O_{\ph(x), Y}\to \Cal
O_{x, X}$ induces the isomorphism of their
completions. If $x\in X_{\sm}$ and $\varphi(x)\in
Y_{\sm}$, then $\varphi$ is \'etale at $x$ iff
$d_x\varphi: T_{x,X}\rightarrow T_{\varphi(x),Y}$
is an isomorphism. If $\ph$ is {\it \'etale } at
each point of $X$, it is called {\it \'etale}.

Let $G$ be a reductive algebraic group.
We refer to \cit!{PV} regarding the
following basic facts.

Let $X$ be an affine algebraic $G$-variety. We
denote by $X/\!\!/G$ the categorical quotient of
$X$ by $G$, i.e\., the affine algebraic variety
whose algebra of regular functions is $k[X]^G$
(it is finitely generated by Hilbert's theorem).
We denote by $\pi^{}_X=\pi^{}_{X,G}:X\to
X/\!\!/G$ the canonical projection (induced by
the inclusion $k[X]^G\hookrightarrow k[X]$).

Every fiber of $\pi^{}_X$ contains a unique
closed orbit. For a point $x\in X$, if the orbit
$G\cdot x$ is closed in $X$, then the stabilizer
$G_x$ is reductive. If $G$ is finite, then
$\pi_X:X\to X/\!\!/G$ is a geometric quotient and
we denote $X/\!\!/G$ by $X/G$.

Let $\ph:X\to Y$ be a $G$-morphism of affine
$G$-varieties. Then we have the following
commutative diagram
$$\CD X @>{\ph}>> Y \\
@V{\pi^{}_{X,G}}VV @VV{\pi^{}_{Y,G}}V\\
X/\!\!/G @>{\ph/\!\!/G}>> Y/\!\!/G,
\endCD\tag{\nmb:{2.1.1}}$$
where $\ph/\!\!/G$ is induced by the restriction
of $\ph^*\!:k[Y]\to k[X]$ to $k[Y]^G$ (if $G$ is
finite, we denote $\ph/\!\!/G$ by $\ph/G$). If
$\ph/\!\!/G$ is \'etale and \thetag{\nmb|{2.1.1}}
is the base change diagram (so that $X=X/\!\!/G
\x^{}_{Y/\!\!/G}Y$), then $\ph$ is called {\it
excellent}. In this case, $\ph$ itself is
\'etale. If $\ph$ is excellent, then for any
$b\in X/\!\!/G$ the restriction of $\ph$ yields a
$G$-isomorphism $ \pi_{X}\i(b)\to
\pi_{Y}\i((\ph/\!\!/G)(b))$. Hence if $Y$ is
irreducible and $\varphi $ is excellent, then $G$
acts faithfully on $X$ iff $G$ acts faithfully on
$Y$.

\subhead\nmb.{2.2.} \'Etale slices
\endsubhead
The following theorem is proved in
\cite{Lu}.

\proclaim{Theorem 2.2.1 } Let $X$ be an affine
algebraic $G$-variety for a reductive algebraic
group $G$, and let $x\!\in\! X$ be a point such
that the orbit $G\!\cdot\! x$ is closed in~$X$.
\roster \item"(i)" There is a $G_x$\!-stable
locally closed affine subvariety $S\!\subseteq\!
X$, called \'etale slice at $x$, such that $x\in
S$ and the $G$-morphism
$$\sigma^{}_S: G*^{}_{G_x}S\to X\tag2.2.1$$
induced by $G\x S\to X$, $(g,s)\mapsto
g\cdot s$, is excellent.
\endroster
\roster \item"(ii)" If $x\in X_{\sm}$, there
exist a $\pi^{}_{S, G_x}$\!-saturated affine open
neighborhood $U$ of $x$ in $S$ and a
$G_x$-equivariant excellent morphism
$$\qquad \la_x:U\to T_{x,S}.\qquad \square \tag2.2.2$$
\endroster
\endproclaim

\vskip -3mm

 Here by
$G*^{}_{G_x}S$ we denote the homogeneous
fiber space over $G/G_x$ with the fiber
$S$, i.e., the categorical quotient of
$G\x S$ by the $G_x$-action
$h\cdot(g,s)=(gh\i,h\cdot s)$ (actually
in this case it is the geometrical
quotient).

Theorem 2.2.1 is crucial for this paper. We will
apply it only in the setting of finite group
actions, in which case its proof is simpler than
in general one. Therefore for making the paper
more self-contained, in short Appendix at the end
of the paper we give, following the referee's
suggestion, some details of the proof of Theorem
2.2.1 for finite $G$.

\subhead\nmb.{2.3}. The Luna
stratification
\endsubhead Let $G$ be a reductive algebraic
group. Denote by $\Cal V(G)$ the set of
isomorphism classes of all algebraic
vector $G$-bundles over homogeneous
spaces $G/H$ where $H$ ranges over all
reductive subgroups of $G$.

Let $X$ be an affine algebraic
$G$-variety. For a point $b\in X/\!\!/G$,
let $G\cdot x$ be the unique closed orbit
in $\pi_X\i(b)$. Since $G_x$ is
reductive, there is a $G_x$-stable direct
linear complement $N_x$ to $T_{x, G\cdot
x}$ in $T_{x, X}$. The $G$-bundle
$G*^{}_{G_x}N_x\to G/G_x$, called the
{\it normal bundle} of $G\cdot x$, is the
representative of some element $\nu
(b)\in \Cal V(G)$. In this way we get a
mapping $\nu:X/\!\!/G \to \Cal V(G)$,
$b\mapsto \nu(b)$.

Let $s\in \Cal V(G)$. Then the sets
$(X/\!\!/G)_s:= \nu\i(s)$ and
$X_s\!:=\!\pi_X\i((X/\!\!/G)_s)$ are
called the {\it strata} of type $s$ of
$X/\!\!/G$ and $X$ respectively.

The following theorem is proved in
\cite{Lu}.

\proclaim{Theorem 2.3.1} Retain the above
notation and assume that $X$ is smooth. Then
\roster \item"(i)" The image of $\nu$ in $\Cal
V(G)$ is finite. \item"(ii)" Each stratum
$(X/\!\!/G)_s$ is a smooth locally closed
subvariety of~$X/\!\!/G$. \item"(iii)" Each
stratum $X_s$ is a locally closed smooth
$\pi^{}_X$-sa\-tu\-rated subvariety of $X$. There
is a $G$-variety $F_s$ such that each fiber
$\pi_{X}^{-1}(b)$, $b\in (X/\!\!/G)_s$, is
$G$-isomorphic to $F_s$ and the restriction
$\pi^{}_X|^{}_{X_s}\!: X_s \to (X/\!\!/G)_s$ is a
fiber bundle with the fiber $F_s$ locally trivial
in \'etale topology {\rm(}i.e., each point of
$(X/\!\!/G)_s$ has an open neighborhood $U$ such
that the pull back of this bundle over a suitable
\'etale covering of $U$ becomes trivial{\rm)}.
\endroster
\endproclaim

If $\ph:X\to Y$ is an excellent morphism
of smooth irreducible affine algebraic
$G$-varieties, then $X_s=\ph\i(Y_s)$ and
$(X/\!\!/G)_s=(\ph/\!\!/G)\i((Y/\!\!/G)_s)$
for every $s\in \Cal V(G)$.

\subhead\nmb.{2.4}. Quotients by finite
group actions
\endsubhead
Now assume that $G$ is a finite group and
$X$ an algebraic $G$-variety (not
necessarily affine). Consider the
property
$$\text{Every $G$-orbit is contained in
an affine open subset of $X$}. \tag Q
$$
For instance, $(Q)$ holds if $X$ is
quasi-projective.

Then one has the following criterion.

\proclaim{Theorem 2.4.1} The following properties
are equivalent:

\roster \item"(i)" $(Q)$ holds. \item"(ii)" There
exists a geometric quotient $$\pi^{}_X=\pi^{}_{X,
G}: X\rightarrow X/G.\tag2.4.1$$
\endroster

The morphism $\pi^{}_X$ is automatically affine
{\rm(}i.e., for any open affine sub\-set $U$ of
$X/G$ the set $\pi^{-1}_X(U)$ is affine{\rm )}.
\endproclaim
\demo{Proof} See \cite{PV, Theorem 4.14}.\quad
$\square$
\enddemo

Assume that there exists a geometric quotient
(2.4.1). Take a point $x\in X$. As the group $G$
is finite, the orbit $G\cdot x$ is closed. As the
morphism $\pi^{}_X$ is affine, the orbit $G\cdot
x$ is contained in a $G$-stable affine open
subset of $X$. Hence Theorem 2.2.1  is applicable
to the action of $G$ on this subset. This yields
\proclaim{Theorem 2.4.2} Assume that for an
action of a finite group $G$ on an algebraic
variety $X$ there exists a geometric quotient
$(2.4.1)$. Then statements {\rm(i)} and {\rm(ii)}
of Theorem $2.2$ hold.
\endproclaim

We retain the notation of Subsection 2.3. For any
point $b\in X/G$, the fiber $\pi_X^{-1}(b)$ is an
orbit $G\cdot x$. As it is finite,
$N_x=T^{}_{x,X}$, so $\nu(b)$ is the class of
$G*^{}_{G_x}T^{}_{x,X}$. Like in Subsection 2.3,
we consider the mapping $\nu: X/G\rightarrow \Cal
V(G)$, $b\mapsto \nu(b)$, and for any $s\in \Cal
V(G)$ define the sets $(X/G)_s:= \nu\i(s)$ and
$X_s\!:=\!\pi_X\i((X/G)_s)$. We call them the
{\it strata} of type $s$ of $X/G$ and $X$
respectively. As every orbit is contained in a
$G$-stable affine open subset of $X$, applying
Theorem 2.3.1  we obtain the following

\proclaim{Theorem 2.4.3} Assume that for an
action of a finite group $G$ on a smooth
algebraic variety $X$ there exists a geometric
quotient $(2.4.1)$. Then statements {\rm(i)},
{\rm(ii)} and {\rm(iii)} of Theorem $2.3$ hold.
\endproclaim

The definition of strata yields the
following description. Let $H$ be the
stabilizer of a point in $X$, and let $M$
be a finite dimensional algebraic
$H$-module over $k$. Then
$$ G\cdot \{\, x\in X^{H}\mid G_x= H, \text{
and the $H$-modules $T_{x, X}$ and $M$
are isomorphic}\, \}$$ is a (possibly
empty) stratum $X_s$, where $s$ is the
class of $G*^{}_{H}M$, and each stratum
of $X$ is obtained in this way. If $X_s$
is nonempty, $Z$ is an irreducible
component of $X_s$, and $z$ is a point of
$Z$, then existence of morphism (2.2.2)
implies that $\dim Z=\dim M^H$ and
$Z\subseteq X^{G_z}$.

We denote by $(X/G)_i$, resp.\,$X_i$, the union
of all strata of $X/G$, resp.\,$X$, of
codimension $i$. This is a locally closed smooth
subvariety of $X/G$, resp.\,$X$. We have
$X_i=\pi^{-1}_X((X/G)_i)$. The subvariety
$(X/G)_0$, resp.\,$X_0$, is a single stratum of
$X/G$, resp.\,$X$. It is called the {\it
principal} stratum and is characterized among all
$(X/G)_i$, resp.\,$X_i$, by the property that it
is dense and open in $X/G$, resp.\,$X$. We have
the inclusion $(X/G)_0\subseteq (X/G)_{\sm}$. If
$X$ is irreducible, $X_0$ coincides with the open
set of all points whose stabilizer is the kernel
of action; in particular if the action is
faithful, $X_0$ is the open set of all points
with trivial stabilizers. By Theorem 2.4.2, the
map $\pi^{}_X|^{}_{X_{\pr}}\!: X_{\pr}\to
(X/G)_{\pr}$ is \'etale.

\subhead\nmb.{2.5}. Local description of
the codimension 1 strata
\endsubhead
The following lemma is a basic tool of
our analysis of pull backs of tensor
fields on quotient varieties. It shows
that, in \'etale topology, a neighborhood
of a point of a codimension 1 stratum of
any quotient variety coincides with such
a neighborhood for a linear action of a
finite cyclic group generated by
pseudo-reflection.

\proclaim{Lemma 2.5.1} Assume that for a faithful
action of a finite group $G$ on an irreducible
smooth algebraic variety $X$ there exists a
geometric quotient $(2.4.1)$. Let $x$ be a point
of  $X_1$. Then: \roster \item"(i)" $G_x$ is a
finite cyclic group. \item"(ii)" There are a
character $\chi\!: G_x\!\!\rightarrow\!
k^{\times}$ of order $|G_x|$ and the local
parameters $u_1,\dots,u_n$ of $X$ at $x$ such
that $u_n$ is a lo\-cal equation of $X_1$ at $x$,
each $u_1,\ldots, u_{n-1}$ is $G_x$-invariant,
and $g\cdot u_n=\chi(g)u_n$ for all $g\in G_x$.
\item"(iii)" $\pi^{}_{X, G}(x)\in (X/G)_{\sm}$.
\item"(iv)" $y:=\pi^{}_{X, G_x}(x)\in
(X/G_x)_{\sm}$. \item"(v)" There are the local
parameters $v_1,\ldots,v_n$ of $X/G_x$ at $y$
such that
$$\pi_{X, G_x}^*(v_1)=u_1,\ \ldots,\ \pi_{X,
G_x}^*(v_{n-1})=u_{n-1},\ \pi_{X,
G_x}^*(v_{n})=u_{n}^{|G_x|}.\tag2.5.1$$
\endroster
\endproclaim

\demo{Proof} By Theorem 2.4.2, there is an
\'etale slice $S$ of $X$ at $x$. As $G$ is
finite, $S$ is a $G_x$-stable neighborhood of
$x$. Hence the set $U$ in (2.2.2) is a
$G_x$-stable neighborhood of $x$ as well. Since
$X$ is irreducible, $U$ is dense in $X$. As the
action of $G$ on $X$ is faithful, this implies
that the action of $G_x$ on $U$ is faithful as
well.

It follows from smoothness of  $X_1$ and $X$ that
$T^{}_{x, X_1}$ is a hyperplane in $T^{}_{x, X}$.
Since the irreducible component of $X_1$
containing $x$ lays in $X^{G_x}$, we have
$T^{}_{x, X_1}\subseteq T_{x, X}^{G_x}$. As the
action of $G_x$ on $U$ is faithful and morphism
(2.2.2) is excellent, the action of $G_x$ on
$T^{}_{x, S}=T^{}_{x, X}$ is faithful as well. As
there is a one-dimensional $G_x$-stable direct
complement to $T^{}_{x, X_1}$ in $T^{}_{x, X}$,
this implies that $T^{}_{x, X_1}=T_{x, X}^{G_x}$,
and $G_x$ is a cyclic group of order $|G_x|$
whose generator acts on $T_{x, X}$ as a
pseudo-reflection with mirror $T_{x, X_1}$.

 From this we deduce that there
are a character $\chi: G_x\rightarrow k^{\times
}$ of order $|G_x|$, and a basis $t_1,\dots,t_n$
of $T_{x, X}^*$ such that $t_1,\dots,t_{n-1}$ are
$G_x$-invariant and $t_n$ is a
$G_x$-semi-invariant with weight $\chi$. The
sequence $t_1,\dots,t_n$ is a system of local
parameters of $T_{x, X}$ at $0$, and $t_n$ is a
local equation of $T_{x, X_1}$ at $0$. Since the
morphism (2.2.2) is \'etale,
$$u_1:=\la_x^*(t_1),
\ldots,u_n:=\la_x^*(t_n) \in k[U]$$ is a
system of local parameters at $x$ that
has the properties stated in (ii).

We have $k[T_{x, X}]=k[t_1,\dots,t_n]$.
Hence $k[T_{x, X}]^{G_x}$ is the
polynomial algebra in $t_1,\dots,t_{n-1},
t_n^{|G_x|}$, and $T_{x, X}/G_x$ is an
affine space.

The set $(G*_{G_x}\!U)/G$ is open in $(G*_{G_x}\!
S)/G$. By Theorem 2.4.2, we have the following
\'etale morphisms
$$ X/G
@<{\sigma_S/G}<<(G\times_{G_x}\!
S)/G\hookleftarrow (G\times_{G_x}\! U)/G=U/G_x
@>{\la_x/G_x}>> T_{x, X}/G_x. $$ Thus the
following commutative diagram arises:
$$
 \CD X
@<{\text{embedding}}<< U @>{\hskip
4.5mm\lambda_x\hskip 4.5mm}>>
T^{}_{x,X}\\
@V\pi^{}_{X\!,\!G}VV
@V\pi^{}_{U\!,\!G_x}VV
@VV\pi^{}_{T_{x\!,\!X}\!,G_x}V\\
X/G @<{\hskip 2.5mm \sigma^{}_S/G\hskip 2.5mm} <<
U/G_x @>{\hskip 1.7mm\lambda_x/G_x\hskip 1.7mm}>>
T^{}_{x,X}/G_x
\endCD
\tag2.5.2
$$

As $\pi^{}_{T_{x\!,\! X}, G_x}(0)\in (T^{}_{x,
X}/G_x)_{\sm}$ and $U/G_x=\pi^{}_{X, G_x}(U)$ is
a neighborhood of $y$ in $X/G_x$, this yields
$y\in (X/G_x)_{\sm}$ and $\pi^{}_{X, G
}(x)=(\sigma_{S}/G)(y)\in (X/G)_{\sm}$.

By our construction, $k[T_{x,
X}/G_x]=k[w_1,\ldots, w_n]$ where $w_1\ldots,w_n$
are defined by the equalities $$\pi_{T_{x,
X}}^*(w_1)=t_1,\ \ldots,\ \pi_{T_{x,
X}}^*(w_{n-1})=t_{n-1}, \ \pi_{T_{x,
X}}^*(w_n)=t_n^{|G_x|}.$$ Since $\lambda_x/G_x$
is \'etale, $v_1:= (\la_x/G_x)^*(w_1),\ldots,
v_n:= (\la_x/G_x)^*(w_n)$ is a system of
parameters of $X/G_x$ at $y$. It follows from
(2.5.2) that
$$\pi_{U, G_x}^*(v_1)=u_1,\ \ldots, \
\pi_{U, G_x }^*(v_{n-1})=u_{n-1},\
\pi_{U, G_x }^*(v_{n})=u_{n}^{|G_x|}.\
\qed $$
\enddemo

\subhead\nmb.{2.6}. Reflection divisor of
quotient variety
\endsubhead Let $X$ be
an irreducible smooth algebraic variety endowed
with a faithful action of a finite group $G$ such
that there exists a geometric quotient (2.4.1).
As $X$ is normal, $X/G$ is normal as well,
\cit!{PV, \S4 and Theorem 3.16}.

It follows from the discussion in Subsection 2.4
that if $(X/G)_1\neq\varnothing$ and
$$(X/G)_1=(X/G)^1_1\cup\ldots\cup
(X/G)^d_1\tag2.6.1$$ is the decomposition
into irreducible components, then for
each $l=1,\ldots, d$ and
$x\in\pi_X\i(b)$, $b\in (X/G)^l_1$ the
integer $|G_x|$ does not depend on $b$
and $x$. Denote it by $r^{}_{l}$.

We encode the information about $r_1,\ldots,r_d$
into the  divisor
$$ R_{X\!/\!G}:=r_{1}R_1+\ldots+r_{d}
R_d\in \Div(X/G),\tag2.6.2 $$ where
$R_l\in\Div(X/G)$ is the prime divisor such that
$\supp(R_l)$ is the closure of $(X/G)^l_1$.

If $(X/G)_1=\varnothing$, we put
$R_{X\!/\!G}=0$. Thus
$$R_{X\!/\!G}=0 \ \Longleftrightarrow\
(X/G)_1=\varnothing. \tag 2.6.3 $$

We call $R_{X\!/\!G}$ the {\it reflection
divisor of} $X/G$. We remark that
$$(X/G)^i_1\cap (X/G)^j_1=\varnothing\ \text{ if }
i\neq j. \tag2.6.4$$
because of smoothness of $(X/G)_1$.

\head\totoc\nmb0{\bf 3}. {\bf Lifting
Tensor Fields} \endhead

\subhead\nmb.{3.1}. Rational tensor fields
on algebraic varieties
\endsubhead Let $Y$ be an irreducible algebraic
variety. Denote by $T(Y_{\sm})$, resp.
$T^*(Y_{\sm})$, the tangent, resp. cotangent,
bundle of $Y_{\sm}$. We define a rational tensor
field of type $\binom {p}{q}$ on $Y$ as a
rational section of the tensor product of $p$
copies of $T(Y_{\sm})$ and $q$ copies of
$T^*(Y_{\sm})$ whose factors are taken in a fixed
order. We call $p+q$ the {\it valency} of such a
section. As changing of the order yields an
isomorphic vector bundle, we consider only the
following representative of this class of vector
bundles:
$$ T^p_q(Y_{\sm}):=T(Y_{\sm})^{\otimes p}\otimes
T^*(Y_{\sm})^{\otimes q}. $$

The set of points where a rational tensor
field $\tau\in
\Gamma_{\operatorname{rat}}(T^p_q
(Y_{\sm}))$ is defined is an open subset
of $Y_{\sm}$. We denote it by $\dom
(\tau)$. We say that $\tau$
is {\it regular} on a subset $U$ of $X$
if $U\subseteq \dom (\tau)$. For a point
$y\in\dom (\tau)$, we denote by $\tau(y)$
the value of $\tau$ at $y$.

The set $\Gamma_{\operatorname{rat}}(T^p_q
(Y_{\sm}))$ of rational tensor fields of type
$\binom {p}{q}$ on $Y$ has a natural structure of
$k(Y)$-module. A rational tensor field $\nu$ of
type $\binom10$ on $Y$ is a rational vector field
on $Y_{\sm}$ or, equivalently, a derivation of
the $k$-algebra $k(Y)$. A rational tensor field
$\omega$ on $Y$ of type $\binom{0}{1}$ is a
rational differential $1$-form on $Y_{\sm}$. We
denote by $\omega(\nu)$ the rational function on
$Y$ whose value at a point $y\in
\dom(\omega)\cap\dom(\nu)$ is equal to the value
of $\omega(y)$ at $\nu(y)$. Rational sections of
$\wedge^{q}(Y_{\sm})$ are rational differential
$q$-forms on $Y_{\sm}$. Rational differential
$q$-forms on $Y$ defined in \cite{Sh} are
naturally identified with them by means of
restriction to $Y_{\sm}$.

Let  $v_1,\dots,v_n$ be a system of local
parameters of $Y_{\sm}$ at a point $y\in
Y_{\sm}$. Then there is a neighborhood
$V$ of $y$ in $Y_{\sm} $ such that
$v_1,\dots,v_n\in k[V]$ and
$dv_1(a),\dots,dv_n(a)$ is a basis of
$T_{a, V}^*$ for each point $a\in V$,
\cite{Sh, Ch. III, \S5}. Let
$\frac{\partial}{\partial v_1},\ldots,
\frac{\partial}{\partial v_n}$ be the
re\-gu\-lar vector fields on $V$ defined
by conditions
$$\textstyle dv_i(\frac{\partial}{\partial
v_j})=\delta_{ij}.\tag\nmb:{3.1.1}$$ Then the
(numbered) set of tensors
$(\frac{\partial}{\partial
v_{i_1}}\otimes\dots\otimes
\frac{\partial}{\partial v_{i_p}}\otimes
dv_{j_1}\otimes\dots\otimes dv_{j_q})(a)$ for
each point $a\in V$ is a basis of $T_{a,
V}^{\otimes p}\otimes T_{a, V}^{* \otimes q}$.

This implies that for each rational tensor field
$\tau\in\Gamma_{\operatorname{rat}}(T^p_q
(Y_{\sm}))$ there is a unique decomposition
$$\textstyle \tau= \sum_{i_1,\dots,i_p,
j_1,\dots,j_q} \tau^{i_1\dots i_p}_{j_1\dots
j_q}\textstyle \frac{\partial}{\partial
v_{i_1}}\otimes\dots\otimes
\frac{\partial}{\partial v_{i_p}}\otimes
dv_{j_1}\otimes\dots\otimes
dv_{j_q},\tag\nmb:{3.1.2}$$ where $\tau^{i_1\dots
i_p}_{j_1\dots j_q}\!\in\! k(Y)$, and we have
$$V\cap\dom (\tau):=\{a\in V\mid \text{ each }
\tau^{i_1\dots i_p}_{j_1\dots j_q} \text{ is
regular at } a \}.\tag\nmb:{3.1.3}$$

>From (3.1.2) we also obtain that for any point
$a\in V$,
$$ \tau(a)=0 \ \Longleftrightarrow\
\tau^{i_1\dots i_p}_{j_1\dots j_q}(a)=0
\text{ for all } i_1,\dots, i_p,
j_1,\dots, j_q. \tag\nmb:{3.1.4} $$ By
(3.1.3), (3.1.4), the set of all zeros of
$\tau$ is a closed subset of $\dom
(\tau)$.

Now assume that the variety $Y$ is
normal. If $\dom (\tau)\neq Y_{\sm}$ and
$$\textstyle Y_{\sm}\setminus \dom
(\tau)=\bigcup_{i=1}^m (\tau)^i_\infty
\tag\nmb:{3.1.5}$$ is the decomposition
of $Y_{\sm}\setminus \dom (\tau)$ into
irreducible components, then by (3.1.3),
 $$\operatorname{codim}^{}_Y (\tau)^i_\infty=
 1 \text{ for
all } i. \tag\nmb:{3.1.6}$$

Thus if the ``variety of poles"
of a tensor field
(i.e., the subvariety of $Y_{\sm}$ where this
field is not defined) is nonempty, it is unmixed
and of codimension 1. On the contrary, if $p$ or
$q\neq 0$, then by (3.1.4), the variety of zeros
of a tensor field of type $\binom{p}{q}$ can be
mixed and/or of codimension greater than $1$.

\subhead\nmb.{3.2}. Divisor of a rational
tensor field
\endsubhead
Let $\tau\in
\Gamma_{\rat}(T_q^p(Y_{\sm}))$ be a
nonzero rational tensor field. We define
$\dvs (\tau)\in \operatorname{Div}(Y)$,
the {\it divisor of} $\tau$, as follows.

Let $C\in\Div(Y)$ be a prime divisor. Since $Y$
is normal, $\supp(C)\cap Y_{\sm}\neq
\varnothing$. Consider in $Y$ an open
neighborhood $U$ of a point $y\in \supp(C)\cap
Y_{\sm}$ where $C$ is given by a local equation
$t$. Then there is a unique integer $\nu(C,
\tau)$ such that \roster \item"(i)" $\supp(C)\cap
\dom (t^{-\nu(C, \tau)}\tau)\neq \varnothing$;
\item"(ii)" $t^{-\nu(C, \tau)}\tau$ does not
identically vanish on $\supp(C)\cap \dom
(t^{-\nu(C, \tau)}\tau)$.
\endroster

To see this, take $V$ and $v_1,\dots,v_n$ as
above assuming moreover that $V\subseteq U$ (this
can be attained by replacing $V$ by $V\cap U$).
Thus (3.1.2) holds. Then (i) and (ii) hold iff
there is a point $z\in \supp(C)$ such that all
functions $t^{-\nu(C, \tau)}\tau^{i_1\ldots
i_p}_{j_1\ldots j_q}$ are defined at $z$ but not
all of them vanish at $z$. This unequi\-vo\-cally
identifies $\nu(C, \tau)$~as $$\nu(C,
\tau)=\min_{i_1,\ldots i_p, j_1,\ldots,j_q} m_{C,
\dvs(\tau^{i_1\ldots i_p}_{j_1\ldots
j_q})}.\tag3.2.1$$

It is immediately seen that $\nu(C,
\tau)$ is well defined (i.e., depends
only on $C$ and $\tau$ but not on the
choice of $V$ and $t$), and that the set
$\{C\mid \nu(C, \tau)\neq 0\}$ is finite
(maybe empty). Now the divisor
$\dvs(\tau)$ is defined by $$m_{C,
\dvs(\tau)}:=\nu(C, \tau)\ \text{ for all
} C.\tag3.2.2$$ For $p=q=0$ this gives
the
usual definition
of divisor of a rational function.

It follows from the definition and
(3.1.5), (3.1.6) that if $m_{C,
\dvs(\tau)}<0$, then $\supp(C)$ is the
closure in $Y$ of one of the
$(\tau)^i_\infty$'s from (3.1.5). Whence
if $\tau$ is a rational tensor field on
$Y$, then
$$\dvs(\tau)\geqslant 0\ \Longleftrightarrow
\tau \text{ is regular on }
Y_{\sm}.\tag\nmb:{3.2.3}$$

On the other hand, we remark that, unlike in the
case of rational functions, in general
$\dvs(\tau)\leqslant 0$ does not imply that
rational tensor field $\tau$ vanishes nowhere.

\subhead\nmb.{3.3}. Pull backs and push forwards
of rational tensor fields
\endsubhead
Consider a dominant rational map $\alpha:
Z\rightarrow Y$ of irreducible algebraic
varieties of the same dimension. As
$\operatorname{char} k=0$, the set $Z_{\alpha}$
of points of $Z_{\sm}$ where $\alpha $ is defined
and \'etale is dense and open in $Z$,
cf.\,\cite{Sh, Ch.\,II, \S 6, no.\,2, Lemma~2}.

Consider a rational tensor field $\tau\in
\Gamma_{\rat}(T^p_q(Y_{\sm}))$. As
$d_z\alpha\!: T_{z, Z}\rightarrow
T_{\alpha(z), Y}$ for $z\in Z_{\alpha}$
is an isomorphism,
 it yields an
isomorphism $$ \alpha^{*}_z :=
((d_z\alpha)^{-1})^{\otimes p} \otimes
((d_z\alpha)^*) ^{\!\otimes q} : T_{\alpha(z),
Y}^{\otimes p}\otimes T_{\alpha(z), Y} ^{*\otimes
q} \rightarrow T_{z, Z} ^{\!\otimes p}\!\otimes
T_{z, Z}^{*\otimes q}.\tag3.3.1$$ Hence for
$W:=\alpha^{-1}(\dom(\tau))\cap Z_{\alpha}$ we
obtain the following section of $T^p_q(W)$:
$$ \gathered \alpha^*(\tau):
W\rightarrow T^p_q(W),\quad
\alpha^*(\tau)(z):=
\alpha^*_z(\tau(\alpha(z))).
\endgathered \tag3.3.2$$

We claim that (3.3.2) is a rational tensor field
on $Z$, regular on $W$. To show this, take a
point $z\in W$, a system of local parameters
$v_1,\ldots, v_n$ of $Y_{\sm}$ at $y=\alpha (z)$,
and a neighborhood $V$ of $y$ in $\alpha(W)$ as
in Subsection 3.1. Then conditions (3.1.1) define
the rational vector fields
$\frac{\partial}{\partial v_j}$ such that
decomposition (3.1.2) holds. Put
$$u_i:=\alpha^*(v_i).\tag3.3.3$$

As $\alpha$ is \'etale et $z$, the sequence
$u_1,\ldots,u_n$ is a system of parameters of
$Z_{\sm}$ at $z$. Let $U$ be a neighborhood of
$z$ in $Z_{\alpha}$ such that $\alpha
(U)\subseteq V$. Then $u_1,\ldots,u_n\in k[U]$
and $du_1(a), \ldots,du_n(a)$ is a basis of
$T^*_{z, U}$ for each point $a\in U$. Let
$\frac{\partial}{\partial
u_1},\ldots,\frac{\partial}{\partial u_n}$ be the
regular vector fields on $U$ defined by the
conditions $du_i(\frac{\partial}{\partial
u_j})=\delta_{ij}$. Then from (3.3.1), (3.3.2)
and (3.1.1) we deduce that
$$\textstyle \alpha^*(dv_i)=du_i \ \text{ and }\
\alpha^*(\frac{\partial}{\partial v_i})=
\frac{\partial}{\partial u_i}\ \text{ for all } \
i.\tag3.3.4$$
%%for all~$i$.

Hence by (3.1.2), there is a unique decomposition
$$\textstyle
\alpha^*(\tau)=\sum_{i_1,\ldots,i_p,j_1,\ldots,j_q}
\alpha^*(\tau^{i_1\ldots i_p}_{j_1\ldots
j_q})\frac{\partial}{\partial
u_{{i_1}}}\otimes
\cdots\otimes\frac{\partial}{\partial
u_{i_p}}\otimes
du_{j_1}\otimes\cdots\otimes
du_{j_q}.\tag3.3.5$$ As $\tau^{i_1\ldots
i_p}_{j_1\ldots j_q}\in k(Y)$, this shows
that $\alpha^*(\tau)$ is a rational
tensor field on $Z$, regular on $W$.
We call $\alpha^*(\tau)$
the {\it pull back} of $\tau$.
Thus we have a map
$$
\alpha^*:
\Gamma_{\rat}(T^p_q(Y_{\sm}))\rightarrow
\Gamma_{\rat}(T^p_q(Z_{\sm})),\quad
\tau\mapsto \alpha^*(\tau).
$$

Now we consider the case where $Z$ is an irreducible smooth algebraic
$G$-va\-riety $X$ for a finite group $G$ such that there exists a geometric
quotient (2.4.1), and $Y=X/G$,
$\alpha=\pi^{}_{X}$. The fixed point set of the
natural action of $G$ on
$\Gamma_{\rat}(T^p_q(X_{\sm}))$ is precisely the
set of all $G$-invariant rational tensor fields
of type $\binom{p}{q}$ on $X$.

If $\tau\in
\Gamma_{\rat}(T^p_q((X/G)_{\sm}))$, then
the definition implies that the field
$\pi_{X}^*(\tau)$ is $G$-invariant. Thus
we obtain a map
$$
\pi_{X}^{*}:
\Gamma_{\rat}(T^p_q((X/G)_{\sm}))
\rightarrow
\Gamma_{\rat}(T^p_q(X_{\sm}))^G, \
\tau\mapsto \pi\!_{X}^{\ *}(\tau).
$$

We claim that $\pi_{X}^{*}$ is bijective. To
prove this, take a field $\vph\!\in\!
\Gamma_{\rat}(T^p_q(X_{\sm}))^G$. As $\dom
(\vph)$ is a $G$-stable and open subset of $X$,
and $X$ is an irreducible variety,
$W:=\pi^{}_X(\dom (\vph))\cap (X/G)_0$ is a dense
open subset of $X/G$. Take points $b\in W$ and
$x\in \pi_X^{-1}(b)$. As $\pi^{}_X$ is \'etale at
each point of $X_0$, we have an isomorphism
$(\pi^{}_X)_{x}^{*}$ (see (3.3.1)). Since $\vph$
is $G$-invariant,
$((\pi^{}_{X})_x^*)^{-1}(\vph(x))$ depends only
on $b$ but not on the choice of $x$ in
$\pi_X^{-1}(b)$. Hence the following formula
defines a section of $T^p_q(W)$:
$$
\pi^{}_{X*}(\vph):
W\rightarrow T^p_q(W),\quad
\pi^{}_{X*}(\vph)(b)=
((\pi^{}_{X})_x^*)^{-1}(\vph(x)).
\tag3.3.6
$$

We claim that $\pi^{}_{X*}(\vph)$ is a rational
tensor field on $X/G$, regular on $W$. We call it
the {\it push forward} of $\vph$. The definitions
imply that the maps $\pi^*_X$ and
$$ \pi^{}_{X*}: \Gamma_{\rat}(T^p_q(X_{\sm}))^G
\rightarrow \Gamma_{\rat}(T^p_q((X/G)_{\sm})), \
\vph\mapsto \pi^{}_{X*}(\vph),
$$
are inverses of one another.

In order to prove the latter claim, fix a system
$v_1,\dots,v_n$ of local parameters of $X/G$ at
$b$, and let $V$ be a neighborhood of $b$ in $W$
such that $v_1,\dots,v_n\in k[V]$ and
$dv_1(b),\ldots, dv_n(b)$ is a basis of
$T_{b,V}^{*}$ for each point $b\in V$. The set
$U:=\pi_X^{-1}(V)$ is an open $G$-stable
neighborhood of $x$ in $\dom (\vph)\cap X_0$. As
$\pi^{}_X|^{}_{X_0}$ is \'etale,
$u_1=\pi^*_X(v_1),\ldots, u_n=\pi^*_X(v_n)\in
k[U]^G$ is  a system of local parameters of $X$
at $x$ such that $du_1(a),\ldots, du_n(a)$ is a
basis of $T_{a,U}^{*}$ for each point $a\in U$.
Hence there is a unique decomposition
$$\textstyle
\vph=\sum_{i_1,\ldots,i_p,j_1,\ldots,j_q}
\varphi^{i_1\ldots i_p}_{j_1\ldots
j_q}\frac{\partial}{\partial u_{{i_1}}}\otimes
\cdots\otimes\frac{\partial}{\partial
u_{i_p}}\otimes du_{j_1}\otimes\cdots\otimes
du_{j_q}, \tag3.3.7$$ where $\varphi^{i_1\ldots
i_p}_{j_1\ldots j_q}\in k(X)$ and
$du_i(\frac{\partial}{\partial
u_j})=\delta_{ij}$.

As the fields $\vph$, $du_i$ and
$\frac{\partial}{\partial u_i}$ are
$G$-invariant, uniqueness of (3.3.7) and the
inclusion $\dom(\vph)\supseteq U$ imply that all
coefficients $\varphi^{i_1\dots i_p}_{j_1\dots
j_q}$ in (3.3.7) are contained in $k[U]^G$.
Hence, as $\pi^{}_X:X\rightarrow X/G$ is a
geometric quotient, there are functions
$\psi^{i_1\dots i_p}_{j_1\dots j_q}\in k[V]$ such
that
$$\varphi^{i_1\dots i_p}_{j_1\dots
j_q}=\pi^*_X(\psi^{i_1\dots i_p}_{j_1\dots
j_q}),\tag{\nmb:{3.3.8}}$$ cf.\,\cit!{PV, \S4}.
On the other hand, the definitions imply that for
any point $a\in U$, $$ \gathered\textstyle
((\pi{}^*_{X})_a)^{-1}\bigl((\frac{\partial}{\partial
u_{i_1}}\otimes\dots\otimes
\frac{\partial}{\partial u_{i_p}}\otimes
du_{j_1}\otimes\dots\otimes du_{j_q})(a)\bigr)
\hskip 2cm
\\ \textstyle
\hskip 3cm
=\bigl(\frac{\partial}{\partial
v_{i_1}}\otimes\dots\otimes
\frac{\partial}{\partial v_{i_p}}\otimes
dv_{j_1}\otimes\dots\otimes
dv_{j_q}\bigr)(\pi^{}_X(a)). \endgathered
\tag{\nmb:{3.3.9}}$$

>From (3.3.6)--(3.3.9), we deduce that the
following decomposition holds $$ \textstyle
\pi_{X*}(\vph)=
\sum_{i_1,\dots,i_p,j_1,\dots,j_q} \psi^{i_1\dots
i_p}_{j_1\dots j_q} \textstyle
\frac{\partial}{\partial
v_{i_1}}\otimes\dots\otimes
\frac{\partial}{\partial v_{i_p}}\otimes
dv_{j_1}\otimes\dots\otimes dv_{j_q}.
\tag{\nmb:{3.3.10}} $$ In turn, (3.3.10) yields
that $\pi_{X*}(\vph)$ is a rational tensor field
of type $\binom{p}{q}$ on $X/G$, regular on $W$.

\subhead\nmb.{3.4}. Regularity of the pull back
of a tensor field on the quotient variety:
Example
\endsubhead
We retain the notation and conventions of the
previous subsection. Our aim in Section 3 is to
give a necessary and sufficient condition of
regularity of the pull back of a rational tensor
field on $X/G$. According to the discussion in
Subsection 3.3, if $\pi\!_{X}^{\ *}(\tau)$ is
regular on $X$, then $\tau$ is automatically
regular on $(X/G)_0$. However in general this
condition is not sufficient: tensor fields on
$X/G$ with ``bad singularities" along the
irreducible components of $(X/G)_1$ may have pull
backs that are not regular on the whole $X$. The
following simplest example illustrates this
phenomenon.

\medskip

\noindent{\bf Example.} Let $X$ be a one
dimensional linear space over $k$. Let
$G$ be a finite group acting on $X$ by
$g\cdot a:=\chi(g)a$ where $\chi:
G\rightarrow k^{\times}$ is a character
of order~$r$. If $u$ is a nonzero linear
function on $X$, then $k[X]=k[u]$ and
$k[X]^G=k[u^r]$. Hence $X/G$ is an affine
line with the coordinate function $v$,
and $\pi^{}_X:=\pi^{}_{X, G}:X\rightarrow
X/G$ is given by $\pi\!_{X}^{\
*}(v)=u^r$.

Each nonzero rational tensor field $\tau$
of type $\binom{p}{q}$ on $X/G$ has the
form
$$\textstyle \tau=\alpha
v^m\prod_{i=1}^{l}(v-\gamma_i)^{d_i}
(\frac{\partial}{\partial v})^{\otimes p}\otimes
(dv)^{\otimes q},$$ where
$m, d_i\in \Bbb Z$, $\alpha, \gamma_i\in
k\setminus \{0\}$
for all $i$, and $\gamma_i\neq \gamma_j$
for all $i\neq j$. Since $\pi\!_X^{\
*}(dv)=ru^{r-1}du$ and $\pi\!_X^{\
*}(\frac{\partial}{\partial v})=
\frac{1}{ru^{r-1}}\frac{\partial}{\partial
u}$, we have $$\gathered\textstyle
\pi\!_X^{\ *}(\tau)=\pi\!_X^{\
*}(v)^m\prod_{i=1}^{l}(\pi\!_X^{\
*}(v)-\gamma_i)^{d_i}(\pi\!_X^{\
*}(\frac{\partial}{\partial v}))^{\otimes
p}\otimes (\pi\!_X^{\ *}(dv))^{\otimes
q}\\= \textstyle
r^{q-p}u^{(r-1)(q-p)+rm}\prod_{i=1}^{l}(u^r-\gamma_i)^{d_i}
(\frac{\partial}{\partial u})^{\otimes
p}\otimes (du)^{\otimes q}.
\endgathered
$$
Therefore the rational tensor field $\pi\!_X^{\
*}(\tau)$ is regular on $X$ if and only if
$$\split \text{ (i) }&\hskip 2mm d_i\geqslant 0
\text{ for all } i,\\ \text{ (ii) }&\hskip 2mm
(r-1)(q-p)+rm\geqslant 0. \endsplit $$

For instance, if $r=3$, $m=-3$, $q-p=7$, and
$d_i=0$ for all $i$, then $\tau$ is regular on
$(X/G)_0$ but not regular on $X/G$, and
$\pi\!_X^{\ *}(\tau)$ is regular on $X$. If
$r=3$, $m=-5$, $q-p=7$ and $d_i=0$ for all $i$,
then $\tau$ is regular on $(X/G)_0$, not regular
on $X/G$, and $\pi\!_X^{\ *}(\tau)$ is not
regular on $X$. \quad $\square $

\medskip

 From the general criterion of regularity of the
pull back that will be obtained below we will see
that this example essentially elucidates the
general mechanism governing the phenomenon under
investigation.

\subhead\nmb.{3.5}. Regularity of the pull back
of a tensor field on a quotient vari\nolinebreak
ety: Case with no strata of codimension $1$
\endsubhead
The corollary of the following
proposition shows that if there are no
strata of codimension $1$, then the
situation  is simple and regularity of
$\tau$ on $(X/G)_0$ is not only necessary
but also sufficient for regularity of
$\pi\!_{X}^{\ *}(\tau)$ on $X$.

\proclaim{Proposition 3.5.1} Let $X$ be an
irreducible smooth algebraic variety endowed with
an action of a finite group $G$. Assume that
there exists a geometric quotient $(2.4.1)$. Let
$U$ be an open subset of $X/G$ such that
$U\cap(X/G)_1=\varnothing$. Consider a tensor
field $\tau\in \Gamma_{\rat}(T^p_q((X/G)_{\sm}))$
such that $U\cap(X/G)_0\subseteq\dom (\tau)$.
Then $\pi_X^{-1}(U)\subseteq\dom (\pi\!_{X}^{\
*}(\tau))$.
\endproclaim

\demo{Proof} As $\dom(\tau)\supseteq
U\cap (X/G)_0$, we have
$\dom(\pi\!_{X}^{\ *}(\tau))\supseteq
\pi_X^{-1}(U)\cap X_0$. Since $U\cap
(X/G)_1=\varnothing$, the complement of
$\pi_X^{-1}(U)\cap X_0$ in
$\pi_X^{-1}(U)$ has codimension
$\geqslant 2$. Now the claim follows from
(3.1.5) and (3.1.6). \qed\enddemo

\proclaim{Corollary 3.5.1} Assume moreover that
$X$ contains no strata of codimension~$1$. Then
the pull back $\pi\!_{X}^{\ *}(\tau)$ of a
rational tensor field $\tau$ on $X/G$ is regular
on $X$ if and only if $\tau$ is regular on the
principal stratum $(X/G)_0$.
\endproclaim

\demo{Proof} As $(X/G)_1=\varnothing$, one can take
$U=X/G$.\qed
\enddemo

\subhead\nmb.{3.6}. Modifying the divisor
of a tensor field by means of a
nonne\-ga\-tive divisor
\endsubhead
Let $Y$ be an irreducible normal
algebraic variety and $B\in \Div(Y)$ a
nonnegative divisor. Consider a nonzero
rational tensor field $\tau\in
\Ga_{\rat}(T^p_q(Y_{\sm}))$. In this
subsection we define a divisor
$\dvs^{}_B(\tau)\in \Div(Y)$ that we call
the {\it $B$-divisor of} $\tau$. It is
obtained from $\dvs (\tau)$ by means of
some ``modification along" $B$. In the
next subsection we will use this notion
in our criterion of regularity of the
pull back of a tensor field $\tau$ on the
quotient variety $X/G$: it is formulated
in terms of
$\dvs^{}_{R_{X\!/\!G}}(\tau)$, where
$R^{}_{X/G}$ is the reflection divisor of
$X/G$, see\,(2.6.2).

If $B=0$, we define
$$\dvs_0(\tau):=\dvs(\tau).\tag3.6.1$$

Now assume that $B>0$,
$$B=b_1B_1+\ldots+b_dB_d,\tag3.6.2$$
where $B_i$ is a prime divisor and $b_i$ a
positive integer for all $i$. Then
$\dvs^{}_B(\tau)$ is obtained from $\dvs(\tau)$
by modifying the multiplicities
$m^{}_{B_1,\dvs(\tau)},
\!\ldots\!,m^{}_{B_d,\dvs(\tau)}$ as defined
below.

To describe this modification, fix some $B_l$. As
$\operatorname{codim}(\supp(B_l))=1$ and $Y$ is
normal, $\supp(B_l)\cap Y_{\sm}\neq\varnothing$.
Hence taking local parameters at a point $y\in
\supp(B_l)\cap \nobreak Y_{\sm}$, we can find an
open neighborhood $V$ of $y$ in $Y_{\sm}$ and a
system $v$ of functions $v_1,\ldots,v_n\in k[V]$
such that $\supp(B_l)\cap \nobreak
V\neq\varnothing$, the ideal of $\supp(B_l)\cap
V$ in $k[V]$ is generated by $v_n$, and
$dv_1(a),\dots,dv_n(a)$ is a basis of $T^*_{a,
V}$ for each $a\in V$. Hence there is a unique
decomposition
$$\textstyle \tau=\sum_{i_1,\dots,i_p,j_1,\dots,j_q}
\textstyle \delta^{i_1\dots
i_p}_{j_1\dots j_q}, \tag{\nmb:{3.6.3}}
$$ where (see (3.1.2))
$$ \textstyle \delta^{i_1\dots
i_p}_{j_1\dots j_q}=\tau^{i_1\dots i_p}_{j_1\dots
j_q}\frac{\partial}{\partial
v_{i_1}}\!\otimes\!\cdots\!\otimes\!
\frac{\partial}{\partial v_{i_p}}\otimes
dv_{j_1}\!\otimes\! \cdots \!\otimes\!
dv_{j_q},\hskip 4mm \tau^{i_1\ldots
i_p}_{j_1\dots j_q}\!\in\!
k(Y).\tag{\nmb:{3.6.4}} $$

With each nonzero tensor field $\delta^{i_1\dots
i_p}_{j_1\dots j_q}$ we associate the following
three integers $m$, $p'$ and $q'$ (depending on
this field):

\medskip
 \roster
\item"$\bullet\hskip 1.3mm     \ $"
$m:=m_{B_{l}, \dvs(\tau^{i_1\ldots
i_p}_{j_1\dots j_q})}$;
\item"$\bullet\hskip 1.3mm     \ $" $p'$
(resp. ~$q'$) is the number of factors
$\frac{\partial}{\partial v_n}$ (resp.
$dv_n$) in the right-hand side of
equality (3.6.4).
\endroster

\medskip

Now define
$$
\gather \mu(b_lB_l, \delta^{i_1\dots
i_p}_{j_1\dots
j_q})^{}_{V,v}:=(b_l-1)(q'-p')+b_lm,
\tag{\nmb:{3.6.5}} \\
\varrho(b_lB_l, \tau)^{}_{V,v} :=\min_{i_1,\dots,
i_p, j_1\dots j_q} \mu(b_lB_l, \delta^{i_1\dots
i_p}_{j_1\dots j_q})^{}_{V,v}. \tag{\nmb:{3.6.6}}
\endgather
$$

It easily follows from (3.6.3)--(3.6.6) that for
any nonzero $f\in k(Y)$ and $\tau,
\varsigma\!\in\! \Gamma_{\rat}(T_q^p(Y)_{\sm})$
such that $\tau+\varsigma\neq 0$ we have
$$\gathered \varrho(b_lB_l, f \tau)^{}_{V,v}
=b_lm^{}_{B_l,
\dvs(f)}+\varrho(b_lB_l, \tau)^{}_{V,v},\\
\varrho(b_lB_l, \tau+\varsigma)^{}_{V,v}\geqslant
\min\{\varrho (b_lB_l, \tau)^{}_{V,v},\, \varrho
(b_lB_l, \varsigma)^{}_{V,v}\}.
\endgathered
\tag{\nmb:{3.6.7}}$$

\proclaim{Lemma 3.6.1} For any nonzero
$\tau\!\in\! \Gamma_{\rat}(T_q^p(Y)_{\sm})$,
$\varsigma\!\in\! \Gamma_{\rat}(T_t^s(Y)_{\sm})$
we~have
$$\varrho(b_lB_l,
\tau\otimes\varsigma)^{}_{V,v}= \varrho(b_lB_l,
\tau)^{}_{V,v}+\varrho(b_lB_l,
\varsigma)^{}_{V,v}.\tag3.6.8$$
\endproclaim
\demo{Proof} If $\tau$ and $\varsigma$ are ``
monomials in $v$", i.e.,
$\tau\!=\!f\!\frac{\partial}{\partial
v_{i_1}}\!\otimes\!\cdots\!\otimes\!
\frac{\partial}{\partial v_{i_p}}\otimes
dv_{j_1}\!\otimes\! \cdots \!\otimes\! dv_{j_q}$,
$\varsigma\!=\!g\!\frac{\partial}{\partial
v_{a_1}}\!\otimes\!\cdots\!\otimes\!
\frac{\partial}{\partial v_{a_s}}\otimes
dv_{b_1}\!\otimes\! \cdots \!\otimes\! dv_{b_t}$
where $f, g\in k(Y)$, then $\tau\otimes
\varsigma$ is such a monomial as well, and the
number $p'$ (resp. $q'$, $m$) for
$\tau\otimes\varsigma$ is the sum of numbers $p'$
(resp. $q'$, $m$) for $\tau$ and $\sigma$.
Therefore in this case the claim readily follows
from the definitions (3.6.5), (3.6.6).

Consider the general case. Then by (3.1.2) we
have $\tau= \sum_i\tau_i$,
$\varsigma=\sum_j\varsigma_j$ where each $\tau_i$
and $\varsigma_j$ is a monomial in $v$. From
(3.6.6) we obtain $\varrho(b_lB_l,\tau)^{}_{V,v}=
\min_i\varrho(b_lB_l,\tau_i)^{}_{V,v}$,
$\varrho(b_lB_l,\varsigma)^{}_{V,v}=
\min_j\varrho(b_lB_l,\varsigma_j)^{}_{V,v}$.
Since $\tau\otimes\varsigma=\sum_{i,
j}\tau_i\otimes\varsigma_j$ and each
$\tau_i\otimes\varsigma_j$ is a monomial in $v$,
we have $$\gathered\varrho(b_lB_l,
\tau\otimes\varrho)^{}_{V,v}= \min_{i,j}
\varrho(b_lB_l,
\tau_i\!\otimes\!\varrho_j)^{}_{V,v}\!=\!
\min_{i,j}\left(\varrho(b_lB_l,
\tau_i)^{}_{V,v}\!+\!\varrho(b_lB_l,
\varsigma_j)^{}_{V,v}\right)\\
= \min_i \varrho(b_lB_l,
\tau_i)^{}_{V,v}+ \min_j\varrho(b_lB_l,
\varrho_j)^{}_{V,v}= \varrho(b_lB_l,
\tau)^{}_{V,v}+ \varrho(b_lB_l,\varrho)^{}_{V,v}.
\quad \square \endgathered
$$
\enddemo

\proclaim{Lemma 3.6.2} $\varrho(b_lB_l,
\tau)^{}_{V, v}$ depends only on $\tau$ and
$b_lB_l$ but not on the choice of $\hskip 1mm
V\hskip -1mm$ and $v_1,\ldots,v_n$.
\endproclaim

\demo{Proof} First we remark that if $V'$ is an
open subset of $V$ and $\supp(B_l)\cap \nobreak
V'\neq \varnothing$, then clearly $\varrho
(b_lB_l, \tau)^{}_{V', v}=\varrho (b_lB_l,
\tau)^{}_{V, v}$.

Now consider another open subset $\tilde V$ of
$Y_{\sm}$, and a system $\tilde v$ of functions
$\tilde v_1,\ldots,\tilde v_n\in k[\tilde V]$
such that $\supp(B_l)\cap \tilde V\neq
\varnothing $, the ideal of $\supp(B_l)\cap
\tilde V$ in $k[\tilde V]$ is generated by
$\tilde v_n$, and $d\tilde v_1(a),\ldots, d\tilde
v_n(a)$ is a basis of $T_{a,\tilde V}^*$ for each
point $a\in \tilde V$. We have to show that
$\varrho (b_lB_l, \tau)^{}_{V,v}=\varrho (b_lB_l,
\tau)^{}_{\tilde V, \tilde v}$. Since
$\supp(B_l)$ is irreducible, $\supp(B_l)\cap
V\cap \tilde V\neq \varnothing$. Hence using the
remark above and replacing $V$ and $\tilde V$ by
$V\cap \tilde V$ we can (and shall) assume that
$V=\tilde V$. This reduces the problem to proving
that
$$\varrho (b_lB_l, \tau)^{}_{V,v}
\leqslant \varrho (b_lB_l, \tau)^{}_{V, \tilde
v}.\tag{\nmb:{3.6.9}}$$

To show that (3.6.9) holds, we need the
divisibility properties of some functions.
Namely, the defining property of $\tilde
v_1,\ldots,\tilde v_n$ implies that there are
unique decompositions
 $$ \textstyle
 dv_i=\sum_{j=1}^ng_{ij}d\tilde v_j, \
\frac{\partial}{\partial
v_i}=\sum_{j=1}^n h_{ij}
\frac{\partial}{\partial\tilde v_j},\quad
\text{ where } \ g_{ij}, h_{ij}\in k[V].
\tag{\nmb:{3.6.10}}$$

Since either of $v_n$ and $\tilde v_n$ is
a generator of the ideal of
$\supp(B_l)\cap V$ in $k[V]$, we have
$$v_n=f\tilde v_n, \text{ where } f
\text{ and } 1/f \in k[V].
\tag{\nmb:{3.6.11}}$$

We claim that $$\tilde v_n  \text{
divides } g_{ni}  \text{ and } h_{in}
\text{ in } k[V] \text{ for all }
i=1,\ldots, n-1. \tag{\nmb:{3.6.12}}$$

To prove this, take a point $a\in
\supp(B_l)\cap V$ and an index $i$,
$1\leqslant i\leqslant n-1$. As $\tilde
v_n(a)=0$, from (3.6.11) we deduce that
$$dv_n(a)= f(a) d\tilde
v_n(a).\tag{\nmb:{3.6.13}}$$

On the other hand, by (3.6.10), we have
$dv_n(a)=\sum_{j=1}^ng_{nj}(a)d\tilde v_j(a)$.
Since $d\tilde v_1(a),\ldots, d\tilde v_n(a)$ is
a basis of $T_{a, V}^*$, the latter equality and
(3.6.13) imply that $g_{ni}$ vanishes on
$\supp(B_l)\cap V$, whence it is divisible by
$\tilde v_n$ in $k[V]$. Further, by (3.1.1), the
value of $dv_n(a)$ on $\frac{\partial}{\partial
v_i}(a)$ is equal to $0$. On the other hand, by
(3.6.13) and (3.6.10), it is equal to the value
of $f(a)d\tilde v_n(a)$ on
$\sum_{j=1}^n\!h_{ij}(a) \frac{\partial}{\partial
{\tilde v}_j}(a)$ that in turn is equal to
$f(a)h_{in}(a)$. By (3.6.11), this implies that
$h_{in}(a)=0$. Thus $h_{in}$ vanishes on
$\supp(B_l)\cap V$, whence it is divisible by
$\tilde v_n$ in $k[V]$. This completes the proof
of (3.6.12).

Now we can proceed to proving (3.6.9). By
%%(3.6.5),
(3.6.6), for each nonzero $\delta_{j_1\ldots
j_q}^{i_1\ldots i_p}$ we have $\varrho (b_lB_l,
\delta_{j_1\ldots j_q}^{i_1\ldots
i_p})^{}_{V,v}=\mu (b_lB_l, \delta_{j_1\ldots
j_q}^{i_1\ldots i_p})^{}_{V,v}$. This and (3.6.7)
imply that it suffices to prove (3.6.9) for the
case when $\tau$ is a monomial in $v$,
$$\textstyle \tau=\frac{\partial}{\partial
v_{i_1}}\!\otimes\!\cdots\!\otimes\!
\frac{\partial}{\partial v_{i_p}}\otimes
dv_{j_1}\!\otimes\! \cdots \!\otimes\!
dv_{j_q}.\tag{\nmb:{3.6.14}}$$

So assume that (3.6.14) holds. Then we prove
(3.6.9) by induction with respect to the valency
$p+q$. First we prove (3.6.9) for $p+q=1$, i.e.,
when $\tau = \frac{\partial}{\partial v_{i}}$ or
$\tau=dv_{i}$ for some $i$. Use the notation of
(3.6.5) for $\delta^{i_1\ldots i_p}_{j_1\ldots
j_q}=\frac{\partial}{\partial v_{i}}$. We have
$m=0$, $q'=0$. If $i<n$, then $p'=0$. If $i=n$,
then $p'=1$. So, by (3.6.5), (3.6.6),
$$
\textstyle \varrho(b_lB_l,
\frac{\partial}{\partial v_{i}})^{}_{V,
v}=\mu(b_lB_l, \frac{\partial}{\partial
v_{i}})^{}_{V, v}=\cases 0 & \text{ if
 }
i<n,\\
1-b_l & \text{ if } i=n.
\endcases
\tag3.6.15
$$

Similarly, using the notation of (3.6.5) for
$\delta^{i_1\ldots i_p}_{j_1\ldots j_q}=dv_i$, we
have: $m=0$, $p'=0$; if $i<n$, then $q'=0$, if
$i=n$, then $q'=1$. By (3.6.5), (3.6.6), this
yields
$$
\varrho(b_lB_l, dv_i)^{}_{V, v}=\mu(b_lB_l,
dv_i)^{}_{V, v}=\cases 0 & \text{ if } i<n,\\
b_l - 1 & \text{ if } i=n.
\endcases \tag 3.6.16
$$

According to (3.6.10), (3.6.6) we have
$$
\gathered
 \varrho(b_lB_l, {\textstyle
\frac{\partial}{\partial v_{i}}})^{}_{V, \tilde
v}=\min_{j} \mu(b_lB_l, h_{ij}{\textstyle
\frac{\partial}{\partial {\tilde v}_{j}}})^{}_{V,
\tilde v},\\
\varrho(b_lB_l, dv_i)^{}_{V, \tilde v}=\min_{j}
\mu(b_lB_l, g_{ij}d\tilde v_j)^{}_{V, \tilde v}.
\endgathered \tag 3.6.17
$$
>From (3.6.7) and (3.6.15), (3.6.16) (where $v_i$,
$v$, $i$ are replaced by $\tilde v_j$, $\tilde
v$, $j$ resp.) we obtain
$$
\eqalign{ &\mu(b_lB_l, h_{ij}{\textstyle
\frac{\partial}{\partial {\tilde v}_{j}}})^{}_{V,
\tilde v}=\cases b_lm_{B_l\!, \dvs(h_{ij})} &
\text{ if
 }
j<n,
\\
b_lm_{B_l\!, \dvs(h_{ij})} +1-b_l & \text{ if }
j=n,
\endcases
\cr &\mu(b_lB_l, g_{ij}d\tilde v_j)^{}_{V, \tilde
v} \hskip .5mm=\cases b_lm_{B_l\!, \dvs(g_{ij})}
&
\text{ if } j<n,\\
b_lm_{B_l\!, \dvs(g_{ij})} +b_l-1 & \text{ if }
j=n.
\endcases
 \cr}
\tag3.6.18
$$
By (3.6.10), (3.6.11) we have $m_{B_l,
\dvs(h_{ij})}\geqslant 0$, $m_{B_l\!,
\dvs(g_{ij})}\geqslant 0$ for all $i, j$, and
$m_{B_l, \dvs(h_{in})}\geqslant 1$, $m_{B_l,
\dvs(g_{ni})}\geqslant 1$ for all $i<n$. Hence
(3.6.18) implies
$$
\eqalign{ & \mu(b_lB_l, h_{ij}{\textstyle
\frac{\partial}{\partial {\tilde v}_{j}}})^{}_{V,
\tilde v}\geqslant\cases 0 & \text{ if
 }
j<n,
\\
1 & \text{ if } j=n, i<n,\\
1-b_l & \text{ if } j=n, i=n,
\endcases
\cr & \mu(b_lB_l, g_{ij}d\tilde v_j)^{}_{V,
\tilde v}\hskip .5mm \geqslant\cases 0 &
\text{ if } j<n, i<n,\\
b_l & \text{ if }
j<n, i=n,\\
b_l-1 & \text{ if } j=n.
\endcases
 \cr}
\tag3.6.19
$$

Since $b_l\geqslant 1$, it follows from (3.6.17),
(3.6.19), (3.6.15), (3.6.16) that
$$
\eqalign{ & \varrho(b_lB_l, {\textstyle
\frac{\partial}{\partial v_{i}}})^{}_{V, \tilde
v}\geqslant \left.\cases 0
& \text{ if } i<n,\\
1-b_l  & \text{ if } i=n,
\endcases \right\}
=\varrho(b_lB_l, {\textstyle
\frac{\partial}{\partial v_{i}}})^{}_{V, v},
 \cr & \varrho(b_lB_l, dv_{i})^{}_{V,
\tilde
v}\hskip .5mm\geqslant \left.\cases 0 & \text{ if } i<n,\\
b_l-1 & \text{ if } i=n.
\endcases \right\}=
\varrho(b_lB_l, dv_i)^{}_{V, v},\cr}
$$
whence (3.6.9) for $p+q=1$.

Now assume that $p+q>1$. Then we can write
$\tau=\tau'\otimes\tau''$ where $\tau'$ and
$\tau''$ are the monomials in $v$ that have
valencies $<p+q$. By Lemma 3.6.1, $$ \gathered
\varrho(b_lB_l,\tau)^{}_{V, v}=
\varrho(b_lB_l,\tau')^{}_{V, v}+
\varrho(b_lB_l,\tau'')^{}_{V, v},\\
\varrho(b_lB_l,\tau)^{}_{V, \tilde v}=
\varrho(b_lB_l,\tau')^{}_{V, \tilde v}+
\varrho(b_lB_l,\tau'')^{}_{V, \tilde v},
\endgathered \tag3.6.20
$$ and by the inductive hypothesis,
$$\gathered
\varrho(b_lB_l,\tau')^{}_{V, v}\leqslant
\varrho(b_lB_l,\tau')^{}_{V, \tilde v},\\
\varrho(b_lB_l,\tau'')^{}_{V, v}\leqslant
\varrho(b_lB_l,\tau'')^{}_{V, \tilde v}.
\endgathered
\tag3.6.21$$

Clearly (3.6.9) follows from (3.6.20), (3.6.21).
\qed\enddemo

\medskip

Given Lemma 3.6.2, we denote $\varrho (b_lB_l,
\tau)_{V,v}$ by $\varrho (b_lB_l, \tau)$.

\smallskip

Now we define the divisor $\dvs^{}_B(\tau)$ for a
positive divisor $B$ given by (3.6.2) as follows.
Let $C\in \Div(Y)$ be a prime divisor, then
$$ m^{}_{C,\dvs^{}_B(\tau)}:=\cases \varrho
(b_lB_l, \tau) \text{ if } C=B_l \text{
for some } l,\\ m^{}_{C, \dvs(\tau)}
\text{ if } C\neq B_l \text{ for all } l.
\endcases
\tag3.6.22$$

Note that (3.6.22), (3.2.1), (3.2.2) and Lemma
3.6.1 imply

\proclaim{Corollary 3.6.1} For any positive
divisor $B$ and nonzero tensor fields
$\tau\!\in\! \Gamma_{\rat}(T_q^p(Y)_{\sm})$,
$\varsigma\!\in\! \Gamma_{\rat}(T_t^s(Y)_{\sm})$
the following equality holds:
$$
\dvs^{}_B(\tau\otimes \varsigma)=
\dvs^{}_B(\tau)+\dvs^{}_B(\varsigma).
$$
\endproclaim

\subhead\nmb.{3.7}.  Behavior at \'etale
morphisms
\endsubhead Let
$\alpha: Z\rightarrow Y$ be a surjective \'etale
morphism of irreducible smooth algebraic
varieties. Let $S$ be a subvariety of codimension
1 in $Y$, and let $f$ be a local equation of $S$
at a point $y\in S_{\sm}$. Then $df(y)\neq 0$,
and as $\alpha$ is \'etale,
$d(\alpha^*(f))(z)\neq 0$ for any point
$z\in\alpha^{-1}(y)$. This implies that
$\alpha^*(f)$ is a local equation at $z$ of the
subvariety $\alpha^{-1}(S)$ of $Z$. Hence if
$C\in \Div (Y)$ and $C'\in \Div(Z)$ are prime
divisors, then $m^{}_{C'\!, \alpha^{*}\!(C)}=0$
or $1$, depending on whether $\supp(C')$ is an
irreducible component of $\alpha^{-1}(\supp(C))$
or not.

 \proclaim{Lemma 3.7.1} Let $\tau$
be a nonzero rational tensor field on $Y$ and let
$B\in \Div(Y)$ be a positive divisor. If $C\in
\Div (Y)$ and $C'\in \Div(Z)$ are prime divisors
such that $m^{}_{C'\!, \alpha^{*}\!(C)}=1$, then
$$ m^{}_{C, \dvs^{}_{\!B}\!(\tau)}=m^{}_{C'\!,
\dvs^{}_{\!\alpha^{\!*}\!(B)}\!
(\!\alpha^{*}(\tau)\!)}. \tag3.7.1
$$
\endproclaim
\demo{Proof} We use the notation of Subsections
3.1, 3.2 and 3.3, so that (3.1.2) and (3.3.5)
hold, and $t$ is a local equation of $C$. Then
$\alpha^*(t)$ is a local equation of $C'$ in an
appropriate neighborhood of a point of
$\supp(C')$. Hence
$$m^{}_{C,\dvs(\tau^{i_1\ldots
i_p}_{j_1\ldots j_q})}=
m^{}_{C',\dvs(\alpha^{*}(\tau^{i_1\ldots
i_p}_{j_1\ldots j_q}))}\tag3.7.2$$ for all
$i_1,\ldots, i_p, j_1,\ldots, j_q$. Then from
(3.1.2), (3.3.4), (3.3.5), (3.2.1), (3.2.2) and
(3.7.2) we deduce that
$$m^{}_{C,\dvs(\tau)}=
m^{}_{C',\dvs(\alpha^{*}(\tau))}.\tag3.7.3$$ From
(3.7.3) and (3.6.22) it follows that (3.7.1)
holds if $C\neq B_i$ for each $B_i$ in
decomposition (3.6.2).

Now assume that $C=B^{}_l$ for some $l$. Then
$b^{}_l=m^{}_{C, B}=m^{}_{C',\alpha^*(B)}$. Use
the notation of Subsection 3.6. By (3.3.3), the
local equations of $C$ and $C'$ in the
appropriate neighborhoods are respectively $v_n$
and $u_n$. Therefore from (3.1.2), (3.3.5),
(3.6.5) and (3.6.6) we deduce the equality
$\varrho(b_lC, \tau)\!=\!\varrho(b_lC',
\alpha^*(\tau))$. By (3.6.22), this yields
(3.7.1). \quad $\square$
\enddemo
\proclaim{Corollary 3.7.1} %%Retain the notation
The following properties are equivalent:
\roster \item"(i)" $
\dvs^{}_{\!B}(\tau)\geqslant 0$,
\item"(ii)"
$\dvs^{}_{\!\alpha^{\!*}\!(B)}
(\alpha^{*}(\tau))\geqslant 0. $
\endroster
\endproclaim

\subhead\nmb.{3.8}. Regularity of the pull back
of a tensor field on a quotient variety: General
criterion
\endsubhead Let $X$ be an
irreducible smooth algebraic variety endowed with
a faithful action of a finite group $G$. Assume
that there exists a geometric quotient (2.4.1).
Recall that we defined by (2.6.2) the reflection
divisor $R_{X/G}$ of $X/G$, $$ R_{X/G}:=r_1
R_1+\ldots+r_d R_d\in \Div(X/G).$$

The following theorem is the main result of
Section 3.

\proclaim{Theorem 3.8.1} Let $X$ be an
irreducible smooth algebraic variety endowed with
an action of a finite group $G$. Assume that
there exists a geometric quotient $(2.4.1)$. Let
$\tau$ be a nonzero rational tensor field on
$X/G$. Then the following properties are
equivalent: \roster
\item"(i)"$\dvs(\pi^{*}_X(\tau))\geqslant 0$,
i.e., the pull back $\pi^{*}_X(\tau)$ of $\tau$
is regular on $X$\!;
\item"(ii)"$\dvs_{R^{}_{X\!/\!G}}(\tau)\geqslant
0$.
\endroster
\endproclaim

\demo{Proof} We can (and shall) assume that the
action is faithful.
By the definitions of
$\dvs(\pi^{*}_X(\tau))$ and
$\dvs_{R^{}_{X\!/\!G}}(\tau)$, either of
conditions (i) and (ii) implies that $\tau$ is
regular on $(X/G)_0$ and $\pi^*_{X}(\tau)$ is
regular on $X_0$. So we can (and shall) assume
that $\tau$ and $\pi^*_{X}(\tau)$ share these
properties.

If $R^{}_{X/G}=0$, the claim follows from (2.6.3)
and Corollary 3.5.1. So we shall assume that
$R^{}_{X/G}>0$, i.e., $(X/G)_1\neq \varnothing$.
As $\dom(\pi^{*}_{X}(\tau))\supseteq X_0$,
regularity of $\pi^{*}_{X}(\tau)$ on $X$ is
equivalent to regularity of $\pi^{*}_{X}(\tau)$
on $X_1$ (see (3.1.5)--(3.1.6)).

Let $x$ be a point of $X_1$ and
$z:=\pi^{}_{X}(x)\in (X/G)_1$. By (2.6.4), there
is a unique irreducible component $(X/G)^l_1$ of
$(X/G)_1$ (see (2.6.1)) containing $z$. As
$\pi^{}_{X}|^{}_W: \pi^{-1}_X(W)\rightarrow W$
for any open subset $W$ of $X/G$ is the
geometrical quotient for the action of $G$ on
$\pi^{-1}_X(W)$, cf.\,\cite{PV, 4.2}, the problem
is local, i.e., we have to prove that
$\pi^{*}_{X}(\tau)$ is regular at $x$ if and only
if $\varrho(r^{}_{l}R^{}_l, \tau)\geqslant 0$,
see (3.6.22), (2.6.2). To that end we apply Lemma
2.5.1.

Let $U$ be the $G_x$-stable neighborhood of $x$
from the proof of this lemma such that the
commutative diagram (2.5.2) holds. Then we have
$$\gathered
\pi^{*}_{X}(\tau)|_U=((\sigma^{}_{S}/G)\circ
\pi^{}_{U,G_x})^{*}
(\tau)=\pi^{*}_{U,G_x}(\theta),\\  \text{ where }
\theta=(\sigma^{}_{S}/G)^{*}(\tau)=
(\pi^{}_{U,G_x})_{*}(\pi^{*}_{X}(\tau)).
\endgathered
\tag3.8.1$$

Let $D\in \Div(U/G^{}_x)$ be the prime divisor
such that $\supp(D)$ is the irreducible component
of $(\sigma^{}_{S}/G)^{-1}((X/G)^{l}_{1})$
containing $y=\pi^{}_{X,G_x}(x)$.  Since
$\sigma^{}_{S}/G$ is \'etale, we have $m^{}_{D,
(\sigma^{}_{S}/G)^{*}(R^{}_{X\!/\!G})}=r^{}_l$,
and by Lemma 3.7.1, our problem is reduced to
proving that the following properties are
equivalent: \roster \item"(a)"
$\pi^{*}_{X}(\tau)|^{}_U$ is regular at $x$,
\item"(b)" $\varrho(r^{}_{l}D, \theta)\geqslant
0$.
\endroster

Let $u_1,\dots, u_n$ and $v_1,\ldots v_n$ be the
local parameters from Lemma 2.5.1. Replacing, if
necessary, $U$ by a smaller neighborhood $U'$ of
$x$, and then $U'$ by $\cap_{g\in G_x}\
g\!\cdot\! U'$, we can (and shall) assume that

\smallskip

\roster \item"$\circ\ $" $U=\pi_{X, G_x}^{-1}(V)$
where $V$ is a neighborhood of $y$ in $X/G_x$,
\item"$\circ\ $" $X_1\cap U$ is irreducible,
\item"$\circ\ $" $u_1,\ldots,u_n\in k[U]$,
\item"$\circ\ $" $v_1,\ldots,v_n\in k[V]$,
\item"$\circ\ $" the ideal of $X_1\cap U$ in
$k[U]$ is generated by $u_n$,
 \item"$\circ\ $"
 the ideal of  $\pi^{}_{U,G_x}(X_1)$ in
 $k[V]$  is generated by $v_n$,
\item"$\circ\ $" $du_1(a),\ldots,du_n(a)$, is a
basis of $T^{*}_{a,U}$ for each $a\in U$,
\item"$\circ\ $" $dv_1(b),\ldots, dv_n(b)$ is a
basis of $T^*_{b, V}$  for each $b\in V$.
\endroster

\smallskip

>From (2.5.1) and (3.1.1) we deduce that
$$\gather\
\pi^{*}_{X, G_x}(dv_i)\!=\!du_i \text{
for } 1\!\leqslant\! i\!\leqslant\! n-1,
\text{ and } \pi_{X,
G_x}^*(dv_n)\!=\!r^{}_lu_n^{r_l-1}du_n,
\tag3.8.2
\\
\textstyle\hskip 1.5mm \pi^{*}_{X,
G_x}(\frac{\partial}{\partial v_i})
     \!=\!\frac{\partial}{\partial u_i}
     \text{ for
} 1\!\leqslant \!i\!\leqslant\! n-1,
\text{ and } \pi^{*}_{X,
G_x}(\frac{\partial}{\partial v_n})\!=\!
\frac1{r_lu_n^{r_l-1}}\frac{\partial}{\partial
u_n}. \tag3.8.3
\endgather
$$

\noindent By Lemma 2.5.1 and
$du_i(\frac{\partial}{\partial
u_j})=\delta_{ij}$, for each $g\in G_x$ we have
$$\gathered g\!\cdot\! du_1=du_1,\ \ldots,\
g\!\cdot\! du_{n-1}=du_{n-1},\ g\!\cdot\! du_n=
\chi(g)du_n,\\ \textstyle
g\!\cdot\!\frac{\partial} {\partial
u_1}=\frac{\partial}{\partial u_1},\ \ \ldots,\ \
g\!\cdot\!\frac{\partial} {\partial
u_{n-1}}=\frac{\partial}{\partial u_{n-1}},\ \
g\!\cdot\!\frac{\partial}{\partial u_n}=
\chi(g)^{-1}\frac{\partial}{\partial u_n}.
\endgathered
\tag{\nmb:{3.8.4}}$$

In view of the properties of
$u_1,\ldots,u_n$, there is a unique
decomposition
$$\textstyle
\pi^*_X(\tau)=\sum_{i_1,\dots,
i_p,j_1,\dots,j_q}\psi^{i_1\dots
i_p}_{j_1\dots j_q}
\frac{\partial}{\partial
u_{i_1}}\otimes\dots\otimes
\frac{\partial}{\partial u_{i_p}}\otimes
du_{j_1}\otimes\dots\otimes
du_{j_q},\tag{\nmb:{3.8.5}} $$ where
$\psi^{i_1\dots i_p}_{j_1\dots j_q}\in
k(U)$.

Since $\pi^*_X(\tau)$ is $G_x$-invariant,
(3.8.4) and the uniqueness of
decomposition (3.8.5) imply that each
summand in the right-hand side of (3.8.2)
is $G_x$-in\-va\-riant.

Consider such a summand $$\textstyle
(\pi^*_X(\tau))^{i_1\ldots i_p}_{j_1\ldots j_q}:=
\psi^{i_1\dots i_p}_{j_1 \dots j_q}
\frac{\partial}{\partial
u_{i_1}}\otimes\dots\otimes
\frac{\partial}{\partial u_{i_p}}\otimes
du_{j_1}\otimes\dots\otimes du_{j_q},
\tag{\nmb:{3.8.6}}$$ and assume that in
\thetag{\nmb|{3.8.6}} there are exactly $p'$
factors $\frac{\partial}{\partial u_n}$ and
exactly $q'$ factors $du_n$. Divide $p'-q'$ by
$r^{}_l:=|G_x|$ (see (2.6.2)) with reminder:
$$
p'-q'=sr_l+t, \ s, t\in \Bbb Z, \ 0\leqslant
t<r_l. \tag3.8.7
$$

 Since the order of $\chi$ is $|G_x|$, and the
field (3.8.6) is $G_x$-invariant, (3.8.4) and
(3.8.7) yield
$$ g\!\cdot\!\psi^{i_1\dots i_p}_{j_1\dots j_q}=
\chi(g)^t\psi^{i_1\dots i_p}_{j_1\dots j_q}\
\text{ for all } g\in G_x . \tag{\nmb:{3.8.8}} $$

\medskip

(a)$\Rightarrow$(b). Assume that
$\pi^{*}_{X}(\tau)|_U$ is regular at $x$.
Shrinking $U$ if necessary we can (and
shall) assume that each $\psi^{i_1\dots
i_p}_{j_1\dots j_q}$ is contained in
$k[U]$. We claim that
$$ \psi^{i_1\dots i_p}_{j_1\dots j_q}=
u_n^t\tilde\psi^{i_1\dots i_p}_{j_1\dots j_q} \
\text{ for some }\ \tilde\psi^{i_1\dots
i_p}_{j_1\dots j_q}\in k[U]^{G_{x}}. \tag3.8.9
$$

Indeed, if $t=0$, then $\psi^{i_1\dots
i_p}_{j_1\dots j_q}\in k[U]^{G_x}$ by (3.8.8), so
(3.8.9) holds. Let $t>0$. As $t<|G_x|$, there is
$g_0\in G_x$ such that $\chi(g_0)^t\neq 1$. Take
a point $a\in X_1\cap U$. Since $a\in X^{G_x}$,
from (3.8.9) we obtain that
$g\!\cdot\!\psi^{i_1\dots i_p}_{j_1\dots j_q}(a)=
\psi^{i_1\dots i_p}_{j_1\dots
j_q}(g^{-1}\!\cdot\!a)= \psi^{i_1\dots
i_p}_{j_1\dots j_q}(a)=\chi(g_0)^t\psi^{i_1\dots
i_p}_{j_1\dots j_q}(a)$, whence $\psi^{i_1\dots
i_p}_{j_1\dots j_q}(a)=0$. Thus $\psi^{i_1\dots
i_p}_{j_1\dots j_q}$ vanishes on $X_1\cap U$.
Hence $\psi^{i_1\dots i_p}_{j_1\dots j_q}=
u_n\tilde\psi^{i_1\dots i_p}_{j_1\dots j_q}$ for
some $\tilde\psi^{i_1\dots i_p}_{j_1\dots j_q}\in
k[U]$. As $k[U]$ is an integral domain, from
this, (3.8.9) and Lemma 2.5.1 we deduce that
$g\!\cdot\!\tilde\psi^{i_1\dots i_p}_{j_1\dots
j_q}=\chi(g)^{t-1}\tilde\psi^{i_1\dots
i_p}_{j_1\dots j_q}$ for all $g\in G_x$. Now the
same arguments can be applied to
$\tilde\psi^{i_1\dots i_p}_{j_1\dots j_q}$, and
so on. Whence the claim.

Now, as $\tilde\psi^{i_1\dots
i_p}_{j_1\dots j_q}\in k[U]^{G_x}$,
we have $\tilde\psi^{i_1\dots
i_p}_{j_1\dots j_q}=\pi_{X,
G_x}^*(\ve^{i_1\dots i_p}_{j_1\dots
j_q})$ for some
$$\ve^{i_1\dots
i_p}_{j_1\dots j_q}\in k[V].\tag3.8.10$$

Hence using (3.8.2), (3.8.3) and (3.8.9), we can
rewrite the field (3.8.6) as follows
$$\textstyle (\pi^*_X(\tau))^{i_1\ldots
i_p}_{j_1\ldots j_q}=\pi_{X, G_
x}^*(r^{p'-q'}_lv_n^{p'-q'-s} \ve^{i_1\dots
i_p}_{j_1\dots j_q} \frac{\partial}{\partial
v_{i_1}}\otimes\dots\otimes
\frac{\partial}{\partial v_{i_p}}\otimes
dv_{j_1}\otimes\cdots\otimes dv_{j_q}).
$$
This equality and (3.8.5), (3.8.6),
(3.8.1) imply that
$$\gathered
\textstyle \theta=\sum_{i_1,\ldots,i_p,
j_1,\ldots
j_q}(\pi^{}_{U,G_x})_{*}(\pi^{*}_X(\tau))^{i_1\ldots
i_p}_{j_1\ldots j_q}\\
\textstyle =\sum_{i_1,\ldots,i_p, j_1,\ldots
j_q}\theta^{i_1\ldots i_p}_{j_1\ldots j_q}
\frac{\partial}{\partial
v_{i_1}}\!\otimes\!\cdots\!\otimes\!
\frac{\partial}{\partial v_{i_p}}\!\otimes
\!dv_{j_1}\!\otimes\!\cdots\!\otimes\! dv_{j_q},
\endgathered
 \tag3.8.11
$$
where
$$
\theta^{i_1\ldots i_p}_{j_1\ldots j_q}=
r^{p'\!-q'}_{l}\! v_n^{p'\!-q'\!-\!s}
\ve^{i_1\dots i_p}_{j_1\dots j_q}. \tag3.8.12
$$

Since $v_n$ is the local equation of $D$ at $y$,
it follows from (3.8.12) and (3.8.10) that
$$
m_{D, \dvs(\theta^{i_1\ldots i_p}_{j_1\dots
j_q})}-p'+q'+s\geqslant 0. \tag3.8.13
$$
Taking into account the identity
$$
\gathered \hskip -3.5cm(r_l-1)(q'-p')+r_lm_{D,
\dvs(\theta^{i_1\ldots i_p}_{j_1\dots j_q})}
\\
\hskip 2cm=r_l\bigl( m_{D, \dvs(\theta^{i_1\ldots
i_p}_{j_1\dots j_q})}-p'+q'+s\bigr)+t,
\endgathered \tag3.8.14
$$
we deduce from (3.8.13), (3.8.7) and $r_l>0$ that
$$
(r_l-1)(q'-p')+r_lm_{D, \dvs(\theta^{i_1\ldots
i_p}_{j_1\dots j_q})}\geqslant 0. \tag3.8.15
$$

Now (3.8.11), (3.8.14) and definitions (3.6.5),
(3.6.6) immediately imply that $\varrho(r^{}_lD,
\theta)\geqslant 0$.

\medskip

(b)$\Rightarrow$(a). Assume that
$\varrho(r^{}_lD, \theta)\geqslant 0$. In view of
the properties of $v_1,\ldots,v_n$, there is a
unique decomposition
$$\gather
\textstyle \theta=
\sum_{i_1,\ldots,i_p,j_1,\ldots,j_q}
\vartheta^{i_1\ldots i_p}_{j_1\ldots j_q},\
 \ \text{where}
\tag3.8.16\\
\textstyle \vartheta^{i_1\ldots i_p}_{j_1\ldots
j_q}=\theta^{i_1\ldots i_p}_{j_1\dots
j_q}\frac{\partial}{\partial
v_{i_1}}\otimes\cdots\otimes
\frac{\partial}{\partial v_{i_p}}\otimes
dv_{j_1}\otimes\cdots\otimes dv_{j_q}
 \tag3.8.17
\endgather
$$
for some $\theta^{i_1\ldots i_p}_{j_1\ldots
j_q}\in k(V)$.

Assume that in (3.8.17) there are exactly $p'$
factors $\frac{\partial}{\partial v_n}$ and
exactly $q'$ factors $dv_n$. Let $s$ and $t$ be
defined by (3.8.7). As $v_n$ is a local equation
of $D$ at $y$,  the definition of
$\varrho(r^{}_lD, \theta)$ (see (3.6.5), (3.6.6))
and the inequality $\varrho(r^{}_lD,
\theta)\geqslant 0$ imply that (3.8.15) holds.
>From here, (3.8.14), (3.8.7) and $r_l>0$ we
deduce that (3.8.13) holds. In turn this implies
that shrinking $V$ if necessary, we can (and
shall) assume that
$$\theta^{i_1\ldots i_p}_{j_1\ldots
j_q}=v^{p'-q'-s}_{n}{\tilde\theta}^{i_1\ldots
i_p}_{j_1\ldots j_q}, \ \text{ where }
{\tilde\lambda}^{i_1\ldots i_p}_{j_1\ldots
j_q}\in k[V].\tag3.8.18$$ Then (3.8.17),
(3.8.18), (3.8.2), (3.8.3), (2.5.1) yield
$$\textstyle \pi^{*}_{U,G_x}
(\vartheta^{i_1\ldots i_p}_{j_1\ldots
j_q})=r^{q'-p'}_{l}u^{t}_n\pi^{*}_{U,G_x}
({\tilde\theta}^{i_1\ldots i_p}_{j_1\ldots j_q})
\frac{\partial}{\partial
u_{i_1}}\otimes\cdots\otimes
\frac{\partial}{\partial u_{i_p}}\otimes
du_{j_1}\otimes\cdots\otimes du_{j_q}.
$$

As $t\geqslant 0$ and ${\tilde\theta}^{i_1\ldots
i_p}_{j_1\ldots j_q}\in k[V]$, this shows that
$\pi^{*}_{U,G_x} (\vartheta^{i_1\ldots
i_p}_{j_1\ldots j_q})$ is regular on $U$. Hence
by (3.8.16) and (3.8.1), the field
$\pi^{*}_{X}(\tau)|^{}_{U}$ is regular on $U$.
\qed\enddemo

\head{\bf 4. Some Applications }\endhead

\subhead\nmb.{4.1. Generalization of
Solomon's theorem} \endsubhead Let $V$ be
a finite dimensional vector space over
$k$ and $G$ a finite subgroup of
$\operatorname{GL}(V)$. It is well known
that the following properties are
equivalent, e.g., see \cite{PV}: \roster
\item "(i)" $V/G$ is smooth. \item "(ii)"
$V/G$ is isomorphic to affine space.
 \item"(iii)" $G$ is generated by
pseudo-reflections.\endroster

Assume that conditions (i)--(iii) hold. Then
Solomon's theorem can be reformulated as the
statement that the algebra of all invariant
regular differential forms on $V$ is the pull
back of the algebra of all regular differential
forms on~$V/G$.

In \cite{M} this statement was generalized to
proper smooth polar actions of Lie groups on
Riemannian manifolds, and in \cite{B} to
algebraic actions of reductive groups on smooth
affine algebraic varieties with smooth
categorical quotients. From Theorem 3.8.1 we
deduce the following

\proclaim{Corollary 4.1.1} Let $X$ be an
irreducible smooth algebraic variety endowed with
an action of a finite group $G$. Assume that
there exists a geometric quotient $(2.4.1)$. Let
$\tau$ be a rational differential form on $X/G$.
Then $\pi^{*}_{X}(\tau)$ is regular on $X$ if and
only if $\tau$ is regular on $(X/G)_{\sm}$.
\endproclaim

\demo{Proof} By Theorem 3.8.1, we have to show
that $\dom(\tau)=(X/G)_{\sm}$ is equivalent to
$\dvs^{}_{R^{}_{X\!/\!G}}(\tau)\geqslant 0$. As
either of these conditions implies
$\dom(\tau)\supseteq(X/G)_0$, we can (and shall)
assume that this inclusion holds. As
$(X/G)_1\subset (X/G)_{\sm}$ by Lemma 2.5.1, it
follows from (3.1.5), (3.1.6) that
$\dom(\tau)=(X/G)_{\sm}$ is equivalent to
$\dom(\tau)\supseteq(X/G)_1$.

We can (and shall) assume that the action is
faithful. We use the notation of Subsection 3.6
for $Y=X$ and $B=R_{X\!/\!G}$. The problem is
local, i.e., we have to show that
$\dom(\tau)\supseteq(X/G)_1^l$ (see (2.6.1)) iff
$\varrho (r^{}_lR^{}_l, \tau)\geqslant 0$. For
each field (3.6.4) we have $p'=p=0$ and
$q'\leqslant 1$ by the skew symmetry condition.
Whence $p'-q'= 0$ or $-1$. In both cases
$p'-q'-s=0$. Hence in our case
$m-p'+q'+s\geqslant 0$ is equivalent to
$m\geqslant 0$, and we are done. \quad $\square $
\enddemo

\subhead\nmb.{4.2. Rational tensor fields
of type $\binom{\otimes p}{\wedge q}$}
\endsubhead Consider a more
general type of rational tensor fields
than differential forms. Namely, let $Y$
be an irreducible algebraic variety.
Consider the vector bundle $
T(Y_{\sm})^{\otimes p} \otimes
\wedge^{q}T^{*}(Y_{\sm})$ over $Y_{\sm}$.
Its rational sections are precisely
rational tensor fields of type $\binom
{p}{q}$ on $Y$ that are skew symmetric
with respect to the covariant entries. We
call them {\it rational tensor fields of
type} $\binom{\otimes p}{\wedge q}$. It
is clear that pull backs and push
forwards of rational tensor fields of
type $\binom{\otimes p}{\wedge q}$ are
rational tensor field of type
$\binom{\otimes p}{\wedge q}$ as well.

\proclaim{Corollary 4.2.1} Let $X$ be an
irreducible smooth algebraic variety endowed with
an action of a finite group $G$. Assume that
there exists a geometric quotient $(2.4.1)$. Let
$\varphi $ be a $G$-invariant regular tensor
field of type $\binom{\otimes p}{\wedge q}$ on
$X$. Then $\pi^{}_{X*}(\varphi)$ is regular on
$(X/G)_{\sm}$.
\endproclaim

\demo{Proof} We have $\varphi=\pi^*_X(\tau)$ for
$\tau:=\pi^{}_{X*}(\varphi)$. Arguing like in the
proof of Corollary 4.1.1 we reduce the problem to
showing that $\varrho(r^{}_lR^{}_l,
\tau)\geqslant 0$ implies $\dom (\tau)\supseteq
(X/G)_1^l$. For each field (3.6.4) we have
$q'\leqslant 1$ by the skew symmetry condition.
Hence $p'-q'\geqslant -1$. If $p'-q'=-1$, resp.
$\geqslant 0$, then $s=-1$, resp. $\geqslant 0$.
Therefore $p'-q'-s\geqslant 0$. Hence
$m-p'+q'+s\geqslant \varrho(r^{}_lR^{}_l,
\tau)\geqslant 0$ implies $m\geqslant 0$, and we
are done. \qed\enddemo

\subhead\nmb.{4.3}. Partially symmetric
tensor fields
\endsubhead Let $Y$ be an irreducible
algebraic variety.
Rational sections of
$\operatorname{S}^{q_1}(T^{*}(Y_{\sm}))\otimes
\cdots\otimes
\operatorname{S}^{q_d}(T^{*}(Y_{\sm}))$
are called {\it rational multi-symmetric
covariant tensor fields of type
$(q_1,\ldots,q_d)$} on $Y$. If $Y$ is
smooth and $\sigma$ is such a field, then
$$
T(\dom (\sigma))^{\oplus d}\rightarrow k,
\ v\mapsto (\sigma(a))(v) \text{ for each
} a\in \dom (\sigma), v\in T_{a,
Y}^{\oplus d},
$$
is a rational function on the algebraic variety
$T(Y)^{\oplus d}$, homogeneous of multi-degree
$(q_1,\ldots,q_d)$ with respect to the natural
diagonal action of $k^{\times}$ on $T(Y)^{\oplus
d}$.  This function is regular on $T(Y)^{\oplus
d}$ iff $\sigma$ is regular on $Y$. Every such
rational function is obtained in this way.
Therefore if $Y$ is smooth and endowed with an
action of a finite group $G$, Theorem 3.8.1
provides a characterization of \linebreak
$G$-in\-variant regular functions on
$T(Y)^{\oplus d}$ in terms of rational
multi-symmetric covariant tensor fields on $Y/G$.

If $Y=V$ is a vector space and $G$ a subgroup of
$\operatorname{GL}(V)$, we have $T(V)^{\oplus
d}=V^{\oplus (d+1)}$ and the action of $G$ on
$T(V)^{\oplus d}$ coincides with the diagonal
action on $V^{\oplus (d+1)}$. Thus Theorem 3.8.1
yields a characte\-ri\-zation of $G$-invariant
polynomials on $V^{\oplus (d+1)}$
%%with respect to the diagonal action of $G$
in terms of rational multi-symmetric
covariant tensor fields on $V/G$. Note
that algebras $k[V\oplus\ldots\oplus
V]^G$ were studied in \cite{Hu} for some
groups $G$ generated by
pseudo-reflections (i.e., when $V/G$ is
smooth).

\head{\bf 5. Lifting Automorphisms of
Quotients}\endhead

\subhead\nmb.{5.1}
\endsubhead Let $X$
be an irreducible algebraic variety endowed with
an action of a finite group $G$ such that there
exists the geometric quotients (2.4.1). If
$\alpha$ is an automorphism of $X/G$ induced by a
$G$-automorphism of $X$, then
$\alpha((X/G)_0)=(X/G)_0$ and
$\alpha_*(R^{}_{X\!/\!G})=\nobreak
R^{}_{X\!/\!G}$.

The following theorem is a conversion of
this statement for linear actions.
%%is given by
%%Note that Theorem 5.1.1
It can be seen as complementing the
results in \cite{B}, \cite{Sch},
\cite{Lo}, and \cite{KLM}.

\proclaim{Theorem 5.1.1} Let $V$ be a finite
dimensional complex vector space and let $G
\subset \operatorname{GL}(V)$ be a finite group.
Let $\psi$ be an automorphism of the algebraic
variety $V/G$ such that $\psi(V/G)_0\subseteq
(V/G)_0$ and $\psi^{}_*
(R^{}_{V\!/\!G})=R^{}_{V\!/\!G}$. Then there is
an automorphism $\varphi$ of the algebraic
variety $V$ such that the following diagram is
commutative
$$
\gathered
\xymatrix{
V\ar[d]_{\pi^{}_{V,G}}\ar[r]^\varphi
&V\ar[d]^{\pi^{}_{V,G}\quad .}\\
V/G\ar[r]^\psi &V/G
}
\endgathered\tag5.1.1
$$
\endproclaim

\demo{Proof} It is proved in \cite{KLM} that
there is an analytic automorphism $\varphi$ of
the analytic space $V$ such that the diagram
(5.1.1) is commutative. Hence the claim
immediately follows from

\proclaim{\nmb.Lemma 5.1.1} Let $X$, $Y$ and $Z$
be complex algebraic varieties. Assume that $Y$
is irreducible and smooth and $X$ is endowed with
an action of a finite group $H$. Let the maps
$\alpha, \beta, \gamma$ in the commutative
diagram
$$ \xymatrix{ &
X\ar[dr]^\alpha &
\\
Y\ar[ur]^\gamma\ar[rr]^\beta && Z}
$$
have the properties: \roster \item"(i)" $\alpha$
and $\beta$ are morphisms, \item"(ii)" every
nonempty fiber of $\alpha $ is an $H$-orbit,
\item"(iii)" $\gamma$ is analytic.
\endroster
Then $\gamma $ is a morphism of algebraic
varieties.
\endproclaim

\demo{Proof} Consider in $X\times
Y$ the subsets
$$\gathered \Delta=\{(x,y)\in X\times Y\mid
\alpha(x)=\beta(y)\},\\
\Gamma=\{(\gamma (y), y)\in X\times
Y \mid y\in Y\}.
\endgathered
$$
Clearly $\Delta$ is Zariski closed. As
$\Gamma$ is the graph of $\gamma$ and $Y$
is smooth, (iii) implies that $\Gamma$ is
an analytic subset of $X\times Y$, cf.\,
\cite{He, Ch.\,IV, no.\,5}; clearly it is
analytically isomorphic to $Y$.

Consider the action of $H$ on $X\times Y$
through the first factor.
By (ii), we have
$$ \textstyle \Delta=\bigcup_{h\in
H}h(\Gamma).\tag5.2.1
$$

Taking into account that the
decompositions of any complex algebraic
variety into irreducible components in
the categories of algebraic varieties and
analytic sets coincide, cf.\,\cite{GR,
Ch.\,V, B, Proposition 1}, we deduce from
irreducibility  of $Y$ that $\Gamma$ is
an irreducible analytic set. Hence
(5.2.1) yields that each $h(\Gamma)$ is
an irreducible component of the analytic
set $\Delta$. Thereby $h(\Gamma)$ is an
irreducible component of the algebraic
variety $\Delta$. Thus $\Gamma$ is a
Zariski closed subset of $X\times Y$.

Let $X@<{\ \pi^{}_{1}}<< X\times
Y@>\pi^{}_{2}>>Y$ be the projections. As
$\pi^{}_{2}|^{}_{\Gamma}: \Gamma \rightarrow Y$
is a bijective morphism of algebraic varieties
and $Y$ is normal, Zariski's Main Theorem implies
that $\pi^{}_{2}|^{}_{\Gamma}$ is an isomorphism
of algebraic varieties. Therefore
$(\pi^{}_{2}|^{}_{\Gamma})^{-1}$ is a morphism
(actually, an isomorphism) of algebraic
varieties. Whence
$\gamma=\pi^{}_{1}\circ(\pi^{}_{2}
|^{}_{\Gamma})^{-1}$ is a morphism as well.
 \qed\enddemo
\enddemo

\subhead\nmb.{5.2}. Remarks \endsubhead By
Lefschetz's principle, see, e.g.,\,\cite{Si,
p.\,164}, Theorem 5.1.1 remains true over any
algebraically closed field of characteristic 0.
Notice that the proof of the result from
\cite{KLM} used in the proof of Theorem 5.1.1
involves lifting of a flat torsion free
connection (that could be carried over to the
algebraic setting) and the geodesic exponential
mapping that in turn involves solving ordinary
differential equations.

For general actions and $k=\Bbb C$ there are
further obstructions for lifting of automorphisms
of $X/G$ that are related to the fundamental
group of $X/G$.

\head{\bf 6. Appendix }\endhead

Here we give some details of the proof of Theorem
2.2.1 for finite $G$.

This proof is based on the following preliminary
results all proved in \cite{Lu}.

\smallskip

(a) Let $X$ be a normal affine variety with a
faithful action of $G$. Let $H$ be a subgroup of
$G$ and let $x$ be a point of $X$. Then the
natural map $\alpha: X/H \to X/G$ is \'etale at
$\pi^{}_{X,H}(x)$ if and only if $H\supseteq
G_x$.

(b) Let in addition $\varphi:X\to Y$ be a finite
equivariant morphism to another normal affine
$G$--variety. If $\varphi$ is \'etale at $x$, and
if $G_{\varphi(x)}=G_x$, then the map
$\varphi/G:X/G\to Y/G$ is \'etale at
$\pi^{}_{X,G}(x)$.

(c) If $x\in X_{\sm}$, then there exists a
$G_x$-equivariant morphism $\varphi:X\to T_{x,
X}$ such that $\varphi$ is \'etale at $x$ and
$\varphi(x)=0$.

\smallskip

Statement (a) is \cite{Lu, Lemme 2}, proved by
references to results of commutative algebra; (b)
is \cite{Lu, Lemme 1}, proved on p.~92, Case (B);
and (c) is \cite{Lu, Lemme, p.~96} (one can
assume that $X$ is a $G$-stable subvariety of a
$G$-module $k^n$ and $x=0$, cf.\,\cite{PV,
Theorem 1.5}; if $\pi: k^n\to T^{}_{x,X}$ is the
projection parallel to a $G_x$-stable complement
to $T^{}_{x,X}$, then clearly one can take
$\varphi=\pi^{}|_X$.)

Note that~(a) implies that $X\to X/G$ is \'etale
whenever $G$ acts freely.

To prove Theorem 2.2.1 (i), note that $X/G_x\to
X/G$ is \'etale at $\pi^{}_{X,G_x}(x)$, by (a).
So there exists a $G_x$-invariant open subset
$S\subseteq X$ such that $S/G_x\to X/G$ is
\'etale. Then its base change $X\times^{}_{X/G}
\!S/G_x\to X$ is \'etale as well. But the
composition of $\psi:G*^{}_{G_x}\!\! S\to
X\times^{}_{X/G}\! S/G_x,~
G_x\!\cdot\!(g,s)\mapsto (g\!\cdot\! s,
G_x\!\cdot\! s)$, with this base change is just
the natural map $G*^{}_{G_x}\! S \to X$, which is
obviously \'etale (since it factors as the open
immersion $G*^{}_{G_x}\! S\to G*^{}_{G_x}\! X$,
followed by the projection $G*^{}_{G_x}\!\! X =
G/G_x \times X\to X$). Thus, $\psi$ is \'etale as
well. On the other hand, the (set-theoretic)
fiber of $\psi$ at $\tilde x=G_x\!\cdot\!(e,
x)\in G*^{}_{G_x}\! S$ is the unique point
$\tilde x$. Thus, shrinking $S$, we may achieve
that $\psi$ is an isomorphism.

To prove Theorem 2.2.1 (ii), consider $\varphi$
from (c). Then $\varphi/G_x: X/G_x\to T_{x,
X}/G_x$ is \'etale at $\pi^{}_{X, G_x}(x)$, by
(b). By the preceding argument, the natural map
$X\to X/G_x\times^{}_{T_{x, X}/G_x}\! \!T_{x, X}$
is \'etale at $x$, and its fiber there is the
unique point $x$. This yields a $G_x$-invariant
open subset $U\subseteq X$ such that $U/G_x\to
T_{x, X}/G_x$ is \'etale, and $U\to
U/G_x\times^{}_{T_{x, X}/G_x}\! \!T_{x, X}$ is an
isomorphism.

\enddocument